\documentclass[11pt]{article}
\usepackage{amssymb}

\newtheorem{Lemma}{Lemma}[section]

\newtheorem{Proposition}[Lemma]{Proposition}

\newenvironment{Proof}
{\begin{trivlist} \item[]{\bf Proof. }}%
{\hspace*{\fill}$\rule{.3\baselineskip}{.35\baselineskip}$\end{trivlist}}

\makeatletter
\@addtoreset{equation}{section}
\makeatother

\newfam\bifam
\font\tenbi=cmmib10 scaled \magstep1
\font\sevenbi=cmmib10 at 11pt
\font\fivebi=cmmib10 at 6pt
\textfont\bifam = \tenbi
\scriptfont\bifam= \sevenbi
\scriptscriptfont\bifam= \fivebi

\usepackage[dvips]{epsfig}
\usepackage{graphicx}

\begin{document}

\title{\bf Discrete vector on-site vortices}

\author{P.G. Kevrekidis$^1$ and D.E. Pelinovsky$^2$ \\
{\small $^{1}$ Department of Mathematics, University of
Massachusetts, Amherst, Massachusetts, 01003-4515, USA} \\
{\small $^{2}$ Department of Mathematics, McMaster University,
Hamilton, Ontario, Canada, L8S 4K1} }

\date{\today}
\maketitle

\begin{abstract}
We study discrete vortices in coupled discrete nonlinear
Schr{\"o}dinger equations. We focus on the vortex cross
configuration that has been experimentally observed in
photorefractive crystals. Stability of the single-component vortex
cross in the anti-continuum limit of small coupling between lattice
nodes is proved. In the vector case, we consider two coupled
configurations of vortex crosses, namely the charge-one vortex in
one component coupled in the other component to either the
charge-one vortex (forming a double-charge vortex) or the
charge-negative-one vortex (forming a, so-called, hidden-charge
vortex). We show that both vortex configurations are stable in the
anti-continuum limit if the parameter for the inter-component
coupling is small and both of them are unstable when the coupling
parameter is large. In the marginal case of the discrete
two-dimensional Manakov system, the double-charge vortex is stable
while the hidden-charge vortex is linearly unstable. Analytical
predictions are corroborated with numerical observations that show
good agreement near the anti-continuum limit but gradually deviate
for larger couplings between the lattice nodes.
\end{abstract}

\section{Introduction}

In the past few years, the developments in the nonlinear optics of
photorefractive materials \cite{fleischer} and of Bose-Einstein
condensates in optical lattices \cite{morsch,morsch2,morsch3} have
stimulated an enormous amount of theoretical, numerical and
experimental activity in the area of discrete nonlinear Hamiltonian
systems. A particular focus in this effort has been drawn to the
prototypical lattice model of the discrete nonlinear Schr{\"o}dinger
(DNLS) equation \cite{ijmpb}. The latter, either as a tight binding
limit \cite{alfimov}, or as a generic discrete nonlinear envelope
wave equation \cite{peyrard} plays a key role in unveiling the
relevant dynamics within the appropriate length and time scales.

One of the principal directions of interest in these lattice systems
consists of the effort to analyze the main features of their
localized solutions. In the particular case of two spatial
dimensions, such structures can be regular discrete solitons
\cite{paperI} or discrete vortices (i.e., structures that have
topological charge over a discrete contour) \cite{paperII}. The
study of these types of coherent structures has made substantial
leaps of progress in the past two years with the numerical and
experimental observation of regular discrete solitons
\cite{esolit,esolit2}, dipole solitons \cite{esolit3},
soliton-trains \cite{esolit4}, soliton-necklaces \cite{esolit6} and
vector solitons \cite{esolit5} in photorefractive crystals and
experimental discovery of robust discrete vortex states
\cite{evort,evort2}, based on earlier theoretical predictions
\cite{malomed,malomed2,malomed3}.

On the other hand, the recent years were marked by the experimental
developments in soft condensed-matter physics of Bose-Einstein
Condensates (BECs). Among the important recent observations one can
single out the experimental illustration of the dark
\cite{bec,bec2,bec3}, bright \cite{bec4,bec5} and gap \cite{bec6}
solitons in quasi-one dimensional BECs. The experimental
capabilities seem to be on the verge of producing similar structures
in a two-dimensional context \cite{ol2d}.

In both of the above contexts (nonlinear optics and atomic physics),
multi-component systems were recently studied due to their relevance
to applications. In particular, the first observations of discrete
vector solitons in nonlinear waveguide arrays were reported in
\cite{christo,christo2}, while numerous experiments with BECs were
directed towards studies of mixtures of different spin states of
$^{87}$Rb \cite{myatt,myatt2} or $^{23}$Na \cite{stamper} and even
ones of different atomic species such as $^{41}$K--$^{87}$Rb
\cite{KRb} and $^{7}$Li--$^{133}$Cs \cite{LiCs}. While the above BEC
experiments did not include the presence of an optical lattice, the
addition of an external optical potential could be manufactured
within the present experimental capabilities \cite{morsch3}.

It is the purpose of the present work to address these recent
features of the physical experiments, namely discrete systems with
multiple components. In particular, we aim at addressing the
fundamental issue of how localized excitations are affected by the
presence of {\em two} components which are coupled (nonlinearly) to
each other. While our results will be presented for the specific
example of two coupled DNLS equations with cubic nonlinearities, we
believe that similar features persist in a variety of other models.
We should note here that rather few studies have focused on the
two-dimensional vector generalization of the DNLS equation
\cite{c7,c7a,c7b}. Among others, we mention the work \cite{molina1}
which was motivated by the experimental system of the nonlinear
waveguide arrays proposed in \cite{christo2}. To the best of our
knowledge, these earlier studies did not address vortices in coupled
discrete systems.

For vortices in coupled systems, a number of interesting questions
emerges concerning the stability of particular vortex configurations
(e.g. the so-called vortex cross \cite{evort,evort2}) including the
case of equal charges in both components and the case of opposite
charges between the two components. The former state has a double
vortex charge, while the latter has a hidden vortex charge. It has
been shown for the continuous NLS equation with cubic-quintic
\cite{desyat} and saturable \cite{fangwei} nonlinearities that these
two states have different stability windows.

In the present setting, we examine the stability of such vortex
structures in the discrete case both analytically and numerically.
We use the method of Lyapunov-Schmidt reductions developed earlier
in \cite{paperII}. This method allows for direct analytical
calculations of eigenvalues of the linear stability problem as
functions of the system parameters (such as the coupling between
adjacent lattice sites and the coupling between the two
components).

Our presentation is structured as follows. In section 2, we
introduce the setup and the vortex cross configurations. In section
3, we study the stability of such configurations in the
one-component model. In section 4, we generalize the vortex cross
configuration to the two-component case and compare our results with
numerical computations of the parameter continuations. In section 5,
we deal with a special Manakov case of the system of two DNLS
equations. Finally, in section 6, we summarize our findings.
Appendix A presents technical details for the case of the
single-component vortex cross.

\section{Setup}

We write the coupled system of discrete nonlinear Schr\"{o}dinger
(DNLS) equations in the form:
\begin{eqnarray}
\label{NLS1} i \dot{u}_{n,m} + \epsilon \left( u_{n+1,m} + u_{n-1,m}
+ u_{n,m+1} + u_{n,m-1} \right) + (|u_{n,m}|^2 + \beta |v_{n,m}|^2 )
u_{n,m} & = & 0, \\
\label{NLS2} i \dot{v}_{n,m} + \epsilon \left( v_{n+1,m} +
v_{n-1,m} + v_{n,m+1} + v_{n,m-1} \right) + (\beta |u_{n,m}|^2 +
|v_{n,m}|^2 ) v_{n,m} & = & 0,
\end{eqnarray}
where $\beta$ is a non-negative parameter for the coupling between
the two components $(u,v)$ and $\epsilon$ is a small non-negative
parameter for the coupling between adjacent lattice sites.
Localized modes of the coupled system (\ref{NLS1})--(\ref{NLS2})
take the form:
\begin{equation}
u_{n,m}(t) = \phi_{n,m} e^{i t}, \qquad v_{n,m}(t) = \psi_{n,m} e^{i
\omega t},
\end{equation}
where $\omega$ is a parameter of time-periodic solutions and
$(\phi_{n,m},\psi_{n,m})$ satisfy the system of nonlinear
difference equations:
\begin{eqnarray}
\label{eq1} (1 - |\phi_{n,m}|^2 - \beta |\psi_{n,m}|^2 ) \phi_{n,m}
= \epsilon \left( \phi_{n+1,m} + \phi_{n-1,m} + \phi_{n,m+1} +
\phi_{n,m-1} \right), \\
\label{eq2} (\omega - \beta |\phi_{n,m}|^2 - |\psi_{n,m}|^2 )
\psi_{n,m} = \epsilon \left( \psi_{n+1,m} + \psi_{n-1,m} +
\psi_{n,m+1} + \psi_{n,m-1} \right).
\end{eqnarray}

\begin{figure}[tbp]
\begin{center}
\epsfxsize=7.0cm \epsffile{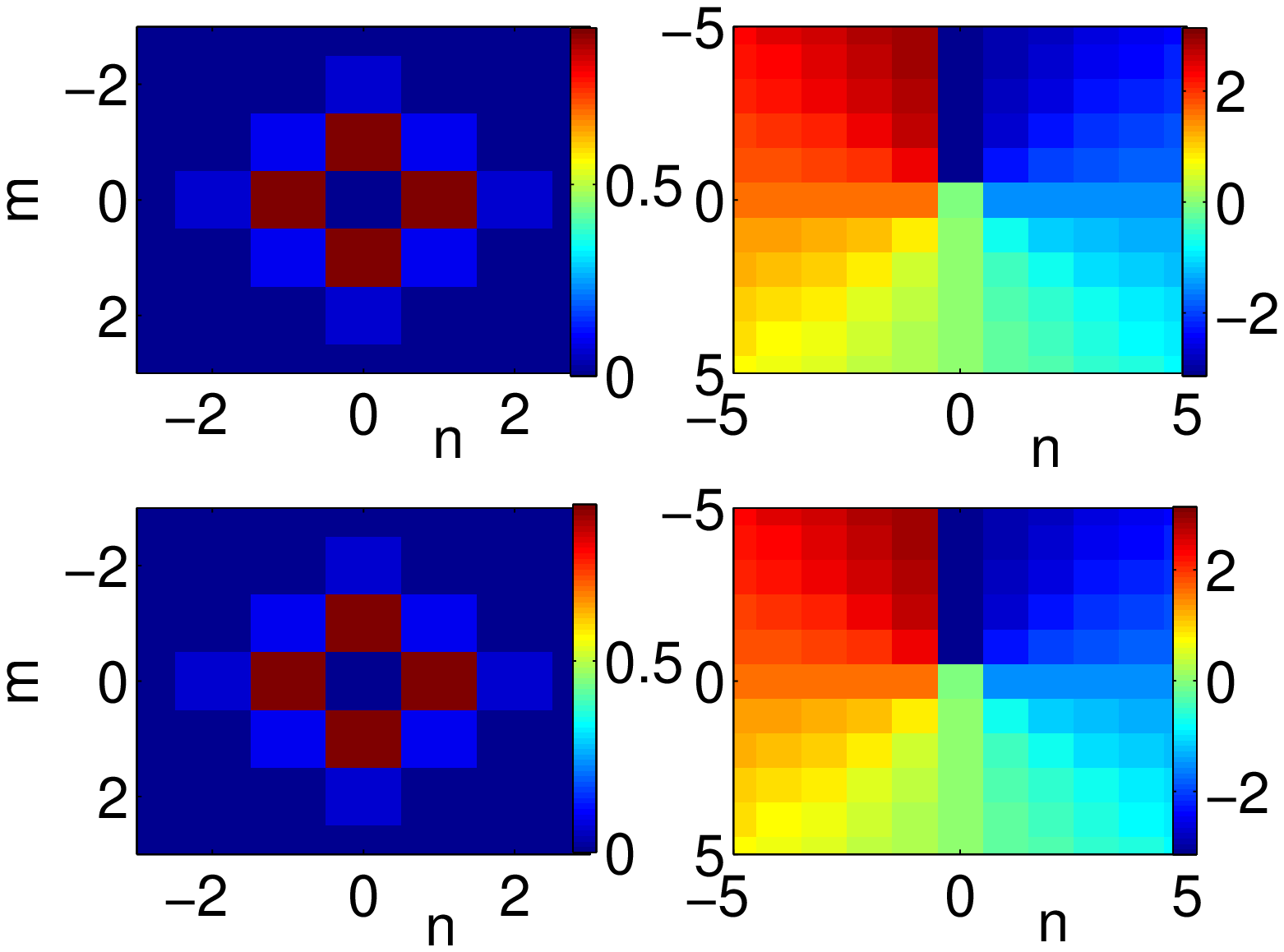} \epsfxsize=7.0cm
\epsffile{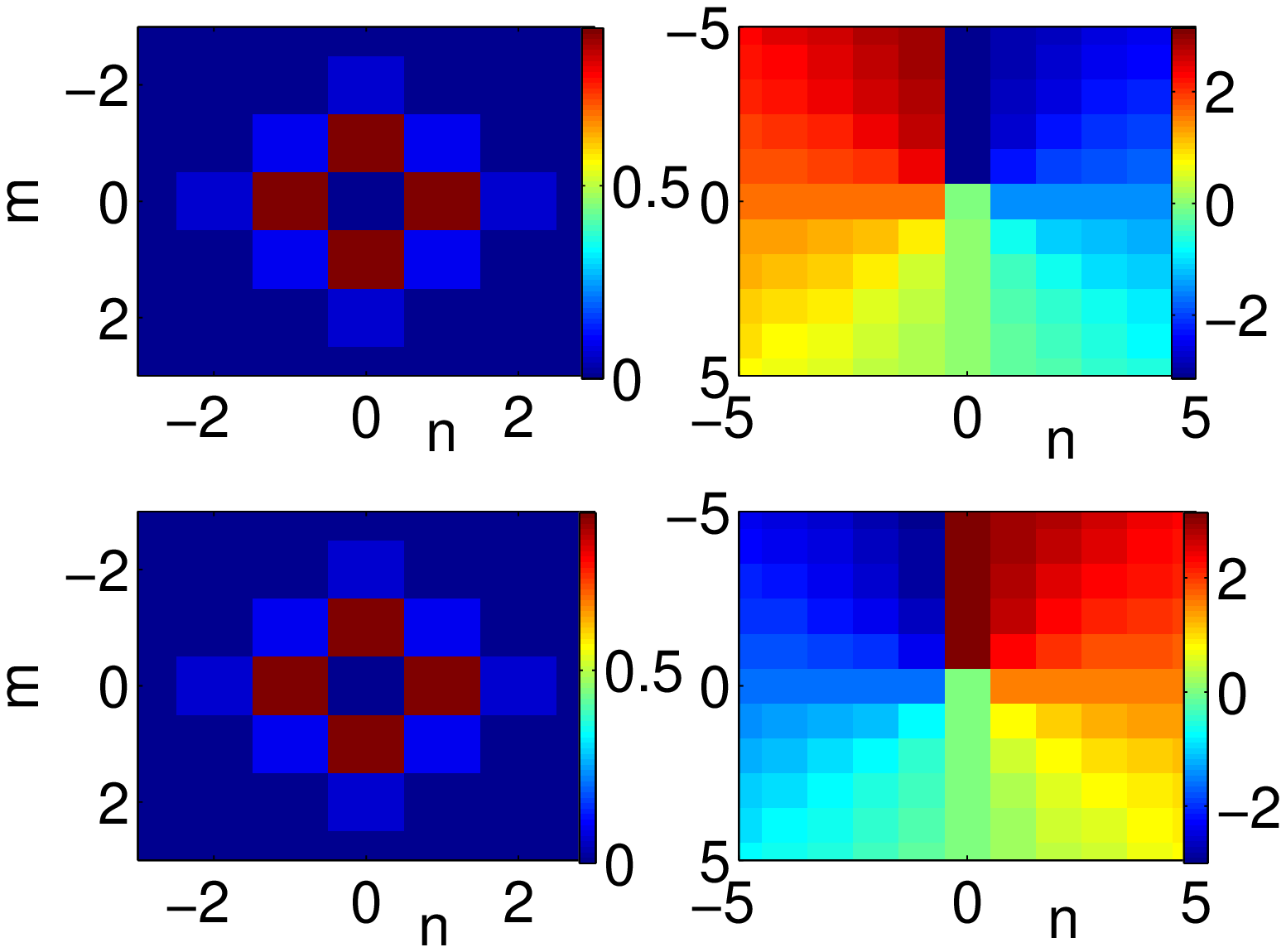} \caption{The contour plots show the
amplitude and phase (left and right panels respectively) of the two
components (top and bottom respectively) for a $(1,1)$ (left four
subplots) and a $(1,-1)$ (right four subplots) vortex configuration,
in the case of $\beta=2/3$, $\omega = 1$, and $\epsilon=0.1$.}
\label{fig0}
\end{center}
\end{figure}

We are interested in a particular vortex solution, called the {\em
vortex cross}. An example of this solution is obtained numerically
for $\beta = \frac{2}{3}$, $\omega = 1$ and $\epsilon = 0.1$ and it
is shown on Figure \ref{fig0}. Let us consider the diagonal square
discrete contour on the grid $(n,m) \in \mathbb{Z}^2$:
\begin{equation}
\label{discrete-contour} S^{(0)} = \{ (-1,0); (0,-1); (1,0); (0,1)
\} \subset \mathbb{Z}^2,
\end{equation}
enumerated in the same order by $j = 1,2,3,4$. We shall assume
that the vortex cross of Figure \ref{fig0} bifurcates from the
limiting solution at the anti-continuum limit $\epsilon = 0$:
\begin{equation}
\label{zero-order-solution} \phi_{n,m}^{(0)} = \left\{
\begin{array}{cc} a e^{i \theta_j}, \quad (n,m) \in S^{(0)} \\ 0,
\quad (n,m) \notin S^{(0)}
\end{array} \right. \qquad \psi_{n,m}^{(0)} = \left\{
\begin{array}{cc} b e^{i \nu_j}, \quad (n,m) \in S^{(0)} \\ 0,
\quad (n,m) \notin S^{(0)}
\end{array} \right.
\end{equation}
where the set of phase parameters $\{\theta_j,\nu_j\}_{j = 1}^4$ is
yet to be determined, while the set of amplitude parameters $(a,b)$
is determined from solutions of the system:
\begin{equation}
\label{quadratic-system} a^2 + \beta b^2 = 1, \qquad \beta a^2 + b^2
= \omega.
\end{equation}
When $\beta \neq 1$, there exists a unique solution of the system
(\ref{quadratic-system}):
\begin{equation}
\label{amplitude-parameters} a^2 = \frac{1 - \beta \omega}{1 -
\beta^2}, \qquad b^2 = \frac{\omega - \beta}{1 - \beta^2}.
\end{equation}
The solution is meaningful only if $a^2 > 0$ and $b^2 > 0$, which
define the domain of existence:
\begin{equation}
\min(\beta,\beta^{-1}) \leq \omega \leq \max(\beta,\beta^{-1}).
\end{equation}
When $\beta = 1$, the domain of existence shrinks into the line
$\omega = 1$ and the solution of the system (\ref{quadratic-system})
forms a one-parameter family:
\begin{equation}
\label{amplitude-parameters-symmetry} a = \cos \delta, \quad b =
\sin \delta, \qquad \delta \in [0,2\pi].
\end{equation}
The vortex cross, if it exists, is defined by the phase
configurations along the discrete contour $S^{(0)}$:
\begin{equation}
\label{symmetric-vortex} \theta_j = \frac{\pi (j-1)}{2}, \qquad
\nu_j = \pm \frac{\pi (j-1)}{2}, \qquad j = 1,2,3,4.
\end{equation}
The upper sign corresponds to the $(1,1)$ coupled state called the
{\em double-charge vortex}, while the lower sign corresponds to
the $(1,-1)$ coupled state called the {\em hidden-charge vortex}.
Persistence and stability of the vortex configurations
(\ref{zero-order-solution}), (\ref{amplitude-parameters}), and
(\ref{symmetric-vortex}) are addressed separately in the cases
$\beta = 0$, $0 < \beta < 1$, $\beta = 1$, and $\beta > 1$.

\section{Scalar vortex cross}

We apply the method of Lyapunov--Schmidt (LS) reductions developed
in \cite{paperII} to the scalar nonlinear difference equation:
\begin{eqnarray}
\label{scalar-eq} (1 - |\Phi_{n,m}|^2 ) \Phi_{n,m} = \epsilon \left(
\Phi_{n+1,m} + \Phi_{n-1,m} + \Phi_{n,m+1} + \Phi_{n,m-1} \right),
\end{eqnarray}
This scalar equation corresponds to the reduction $\psi_{n,m} = 0$,
$\forall (n,m) \in \mathbb{Z}^2$ of the system
(\ref{eq1})--(\ref{eq2}). Local existence of a single-component
vortex cross in the scalar problem (\ref{scalar-eq}) is proved in
Appendix A for small values of $\epsilon$ (on the basis of
Proposition 2.9 in \cite{paperII}). This result is formulated as
follows.

\begin{Proposition}
\label{proposition-existence} There exists a unique (up to the gauge
invariance) continuation in $\epsilon$ of the limiting solution at
$\epsilon = 0$:
\begin{equation}
\label{scalar-vortex-0} \Phi_{n,m}^{(0)} = \left\{ \begin{array}{ll}
e^{i \theta_j}, \quad (n,m) \in S^{(0)} \\ 0, \quad (n,m) \notin
S^{(0)}
\end{array} \right.
\end{equation}
where $S^{(0)}$ is given by (\ref{discrete-contour}) and the
values of $\theta_j$ are given by (\ref{symmetric-vortex}). The
family of vortex solutions $\Phi_{n,m}(\epsilon)$, $(n,m) \in
\mathbb{Z}^2$ is a smooth (real analytic) function of $\epsilon$.
\end{Proposition}

To address spectral stability of the vortex cross in the
time-evolution of the single-component DNLS equation, we consider
the linearization problem with the explicit formula
$$
u_{n,m}(t) = e^{i t} \left[ \Phi_{n,m} + a_{n,m} e^{\lambda t} +
\bar{b}_{n,m} e^{\bar{\lambda} t} \right],
$$
and derive the linear eigenvalue problem from the DNLS equation,
\begin{eqnarray}
\label{stability-1} (1 - 2 |\Phi_{n,m}|^2) a_{n,m} - \Phi_{n,m}^2
b_{n,m} - \epsilon (a_{n+1,m} + a_{n-1,m} + a_{n,m+1} + a_{n,m-1} )
= i \lambda a_{n,m} \\
\label{stability-2} (1 - 2 |\Phi_{n,m}|^2) b_{n,m} -
\bar{\Phi}_{n,m}^2 a_{n,m} - \epsilon (b_{n+1,m} + b_{n-1,m} +
b_{n,m+1} + b_{n,m-1} ) = - i \lambda b_{n,m},
\end{eqnarray}
where $\lambda$ is an eigenvalue and $(a_{n,m},b_{n,m})$ are
components of an eigenvector. Symbolically, we write the linear
eigenvalue problem as
\begin{equation}
{\cal H}(\epsilon) \mbox{\boldmath $\varphi$} = i \lambda \sigma
\mbox{\boldmath $\varphi$}, \label{scalar-stability-problem}
\end{equation}
where ${\cal H}(\epsilon)$ is the linearized Jacobian matrix for the
system (\ref{scalar-eq}), $\sigma$ is a diagonal matrix of $(1,-1)$,
and $\mbox{\boldmath $\varphi$}$ is an eigenvector consisting of
$(a_{n,m},b_{n,m})$. The linear eigenvalue problem for the limiting
solution $\Phi_{n,m} = \Phi_{n,m}^{(0)}$ at $\epsilon = 0$ has a set
of double zero eigenvalues with the eigenvectors ${\bf e}_j$ and
generalized eigenvectors $\hat{\bf e}_j$, such that ${\cal H}^{(0)}
{\bf e}_j = {\bf 0}$ and ${\cal H}^{(0)} \hat{\bf e}_j = 2 i \sigma
{\bf e}_j$, where ${\cal H}^{(0)} = {\cal H}(0)$. The index $j$
enumerates the set $S^{(0)}$ and the eigenvectors ${\bf e}_j$ and
$\hat{\bf e}_j$ have non-zero components only at the corresponding
nodes of the set $S^{(0)}$,
\begin{equation}
{\bf e}_j = i \left( \begin{array}{cc} e^{i \theta_j} \\ - e^{-i
\theta_j} \end{array} \right), \qquad \hat{\bf e}_j = \left(
\begin{array}{cc} e^{i \theta_j} \\ e^{-i \theta_j} \end{array}
\right).
\end{equation}
The kernel of ${\cal H}(\epsilon)$ for $\epsilon \neq 0$ includes at
least one eigenfunction
\begin{equation}
\label{geometric-kernel} \mbox{\boldmath $\varphi$}_{n,m} = \left(
\begin{array}{cc} \Phi_{n,m} \\ - \bar{\Phi}_{n,m} \end{array}
\right), \qquad (n,m) \in \mathbb{Z}^2,
\end{equation}
which follows from the gauge invariance of the DNLS equation with
respect to rotation of the complex phase in $\Phi_{n,m}$, $(n,m) \in
\mathbb{Z}^2$. It is easy to show that a generalized kernel for zero
eigenvalue is non-empty as it includes a solution of the
inhomogeneous equation ${\cal H}(\epsilon) \tilde{\mbox{\boldmath
$\varphi$}} = 2 i \sigma \mbox{\boldmath $\varphi$}$ exists, where
$\mbox{\boldmath $\varphi$}$ is given by (\ref{geometric-kernel}).

Using the perturbation series expansion for $\Phi_{n,m}(\epsilon)$,
we define the expansion ${\cal H}(\epsilon) = {\cal H}^{(0)} +
\epsilon {\cal H}^{(1)} + \epsilon^2 {\cal H}^{(2)} + {\rm
O}(\epsilon^3)$. By Lemma 4.1 in \cite{paperII}, computations of
Appendix A determine the splitting of zero eigenvalues of ${\cal
H}(\epsilon)$ as $\epsilon \neq 0$. The splitting of zero
eigenvalues of $\sigma {\cal H}(\epsilon)$ is formulated and proved
as follows.

\begin{Proposition}
\label{proposition-stability} Let $\Phi_{n,m}(\epsilon)$, $(n,m) \in
\mathbb{Z}^2$ be a family of vortex solutions defined by Proposition
\ref{proposition-existence}. The linearized problem
(\ref{stability-1})--(\ref{stability-2}) has zero eigenvalue of
algebraic multiplicity {\em two} and geometric multiplicity {\em
one} and {\em three} small pairs of purely imaginary eigenvalues of
negative Krein signatures\footnote{A simple eigenvalue of the linear
eigenvalue problem (\ref{scalar-stability-problem}) is said to have
the negative Krein signature if the quadratic form for the
associated eigenvector $({\cal H}(\epsilon) \mbox{\boldmath
$\varphi$},\mbox{\boldmath $\varphi$})$ is negative.} with the
asymptotic approximations,
$$
\lambda_{1,2}, \lambda_{3,4} = \pm 2 i \epsilon + {\rm
O}(\epsilon^2), \qquad \lambda_{5,6} = \pm 4 i \epsilon^2 + {\rm
O}(\epsilon^3).
$$
The rest of the spectrum is bounded away the origin as $\epsilon
\to 0$ and it is located on the imaginary axis of $\lambda$.
\end{Proposition}

\begin{Proof}
We supplement the general proof of Lemma 4.2 in \cite{paperII}
with the explicit perturbation series expansions for small
eigenvalues of the linear eigenvalue problem
(\ref{scalar-stability-problem}):
\begin{equation}
\label{pertubation-eigenvectors-eigenvalues} \mbox{\boldmath
$\varphi$} = \mbox{\boldmath $\varphi$}^{(0)} + \epsilon
\mbox{\boldmath $\varphi$}^{(1)} + \epsilon^2 \mbox{\boldmath
$\varphi$}^{(2)} + {\rm O}(\epsilon^3), \qquad \lambda = \epsilon
\lambda_1 + \epsilon^2 \lambda_2 + {\rm O}(\epsilon^3),
\end{equation}
where
$$
\mbox{\boldmath $\varphi$}^{(0)} = \sum_{j = 1}^4 c_j {\bf e}_j,
\qquad \mbox{\boldmath $\varphi$}^{(1)} = \frac{\lambda_1}{2}
\sum_{j = 1}^4 c_j \hat{\bf e}_j + \mbox{\boldmath
$\varphi$}^{(1)}_{\rm inhom},
$$
and the solution $\mbox{\boldmath $\varphi$}^{(1)}_{\rm inhom} = -
{\cal H}^{(0) -1} {\cal H}^{(1)} \mbox{\boldmath $\varphi$}^{(0)} =
- {\cal H}^{(1)} \mbox{\boldmath $\varphi$}^{(0)}$ is uniquely
defined on the set $S^{(1)}$, where $S^{(1)}$ is the set of adjacent
nodes to the set of $S^{(0)}$. At the second-order perturbation
theory, the problem is written in the form,
\begin{equation}
\label{second-order-pertubation} {\cal H}^{(0)} \mbox{\boldmath
$\varphi$}^{(2)} + {\cal H}^{(1)} \mbox{\boldmath $\varphi$}^{(1)} +
{\cal H}^{(2)} \mbox{\boldmath $\varphi$}^{(0)} = i \lambda_1 \sigma
\mbox{\boldmath $\varphi$}^{(1)} + i \lambda_2 \sigma
\mbox{\boldmath $\varphi$}^{(0)}.
\end{equation}
Projecting the problem to the kernel of ${\cal H}^{(0)}$, we find
the reduced eigenvalue problem:
\begin{equation}
{\cal M}_2 {\bf c} = \frac{1}{2} \lambda_1^2 {\bf c},
\end{equation}
where ${\bf c} = (c_1,c_2,c_3,c_4)^T$ and ${\cal M}_2$ is computed
in Appendix A. Therefore, two negative eigenvalues $\gamma_{1,2}$ of
the Jacobian matrix $\epsilon^2 {\cal M}_2$ generate two pairs of
imaginary eigenvalues of negative Krein signatures in the linear
eigenvalue problem by virtue of the relation $\lambda = \pm \sqrt{2
\gamma}$. The same computation is then extended up to the fourth
order, where it is found that the negative eigenvalue $\gamma_3$ of
the extended matrix $\epsilon^2 {\cal M}_2 + \epsilon^4 {\cal M}_4$
determines the third pair of purely imaginary eigenvalues by virtue
of the same relation $\lambda = \pm \sqrt{2 \gamma}$.
\end{Proof}

We note that the count of eigenvalues of negative Krein signatures
corresponds to the closure theorem for negative index of ${\cal
H}(\epsilon)$ (see \cite{paperI} for details). There are four
negative eigenvalues of ${\cal H}^{(0)}$ for the limiting solution
(\ref{scalar-vortex-0}) and three more small negative eigenvalues
occur for $\epsilon \neq 0$. The total number of negative
eigenvalues is reduced by the gauge symmetry constraint, such that
six negative eigenvalues in a constrained subspace match three pairs
of imaginary eigenvalues with negative Krein signature.

\begin{figure}[tbp]
\begin{center}
\epsfxsize=7.0cm \epsffile{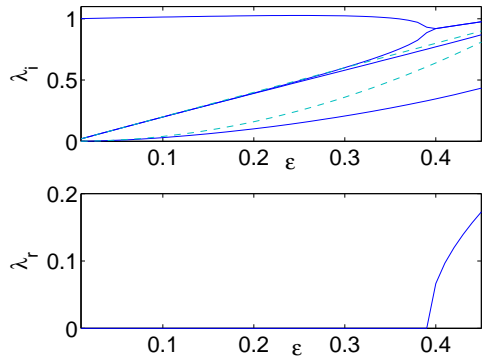} \caption{Eigenvalues of
the scalar vortex cross versus $\epsilon$. The top panel shows the
imaginary part of the relevant eigenvalues, while the bottom panel
shows the real part. The solid lines display the numerical results,
while the dashed ones correspond to the asymptotic approximations.}
\label{fig1}
\end{center}
\end{figure}

The asymptotic approximations of eigenvalues $\lambda$ are plotted
on Figure \ref{fig1} by dashed lines. The numerical computations of
the same eigenvalues (up to the prescribed numerical accuracy)
versus $\epsilon$ are shown by solid lines. All three pairs of
purely imaginary eigenvalues bifurcate into complex domain when they
collide to other eigenvalues of stability problem (e.g. with
eigenvalues of positive Krein signatures or with the spectral band).
The first collision is numerically detected to occur at $\epsilon
\approx 0.395$.

\section{Vector vortex crosses for $0 < \beta < 1$ and $\beta > 1$}

In order to consider the coupled vortex configurations in the
non-degenerate case $\beta \neq 1$, we extend computations of
Appendix A to the solution of the coupled nonlinear difference
equations (\ref{eq1})--(\ref{eq2}). We report here computations for
two related problems: (i) bifurcations of small eigenvalues of the
linearized Jacobian matrix near the zero eigenvalue and (ii)
bifurcations of small eigenvalues of the linearized stability
problem near the origin. Because of the computational complexity of
the analytical approximations, we shall complement the analytical
results of the second-order Lyapunov--Schmidt (LS) reductions with
the symbolic computational results of the fourth-order LS
reductions.

Similarly to the scalar case, the linearized stability problem for
the two-component system takes the matrix-vector form:
\begin{equation}
\label{stability-vector} {\cal H}(\epsilon) \mbox{\boldmath
$\varphi$} = i \lambda \sigma \mbox{\boldmath $\varphi$},
\end{equation}
where ${\cal H}(\epsilon)$ is the linearized Jacobian matrix for the
system (\ref{eq1})--(\ref{eq2}), $\sigma$ is a diagonal matrix of
$(1,-1,1,-1)$, and $\mbox{\boldmath $\varphi$}$ is an eigenvector
consisting of four elements of the perturbation vector at each node
$(n,m) \in \mathbb{Z}^2$. The diagonal block of the matrix ${\cal
H}(\epsilon)$ at each node $(n,m) \in \mathbb{Z}^2$ takes the form:
{\small
$$
{\small  \left( \begin{array}{cccc} 1 - 2 |\phi_{n,m}|^2 - \beta
|\psi_{n,m}|^2 & - \phi_{n,m}^2 & - \beta \phi_{n,m}
\bar{\psi}_{n,m} & - \beta \phi_{n,m} \psi_{n,m} \\
- \bar{\phi}_{n,m}^2 & 1 - 2 |\phi_{n,m}|^2 - \beta |\psi_{n,m}|^2 &
- \beta \bar{\phi}_{n,m} \bar{\psi}_{n,m} & - \beta \bar{\phi}_{n,m} \psi_{n,m} \\
- \beta \bar{\phi}_{n,m} \psi_{n,m} & - \beta \phi_{n,m} \psi_{n,m}
&
\omega - \beta |\phi_{n,m}|^2 - 2 |\psi_{n,m}|^2 & - \psi_{n,m}^2 \\
- \beta \bar{\phi}_{n,m} \bar{\psi}_{n,m} & - \beta \phi_{n,m}
\bar{\psi}_{n,m} & - \bar{\psi}_{n,m}^2 & \omega - \beta
|\phi_{n,m}|^2 - 2 |\psi_{n,m}|^2 \end{array} \right)}
$$
}The non-diagonal blocks of ${\cal H}(\epsilon)$ comes from the
difference operators in the right-hand-side of the system
(\ref{eq1})--(\ref{eq2}).

\subsection{Bifurcations of zero eigenvalues of the linearized Jacobian
matrix}

We extend the perturbation series expansions
(\ref{appendix-expansions}) to the two-component case,
\begin{equation}
\label{vector-expansions} \phi_{n,m}(\epsilon) = \sum_{k =
0}^{\infty} \epsilon^k \phi_{n,m}^{(k)}, \qquad \psi_{n,m}(\epsilon)
= \sum_{k = 0}^{\infty} \epsilon^k \psi_{n,m}^{(k)},
\end{equation}
where the zero-order solution in the anti-continuum limit is given
by (\ref{zero-order-solution}) and parameters $(a,b)$ are given in
(\ref{amplitude-parameters}). The first-order corrections are found
from the uncoupled system of equations, similarly to the scalar
case:
\begin{equation}
\label{first-set-vector} \phi_{n,m}^{(1)} = \left\{
\begin{array}{ll} 0, \quad (n,m) \in S^{(0)}
\\ a \sum_l^{(1)} e^{i \theta_l}, \quad (n,m) \in S^{(1)} \\ 0, \quad (n,m)
\notin S^{(0)} \cup S^{(1)} \end{array} \right. \qquad
\psi_{n,m}^{(1)} =  \left\{ \begin{array}{ll} 0, \quad (n,m) \in
S^{(0)} \\ \omega^{-1} b  \sum_l^{(1)} e^{i \nu_l}, \quad (n,m) \in S^{(1)} \\
0, \quad (n,m) \notin S^{(0)} \cup S^{(1)}
\end{array} \right.
\end{equation}
where $\sum_l^{(1)}$ is defined in (\ref{first-set}). The
second-order corrections are found in the form:
\begin{equation}
\label{second-set-vector} \phi_{n,m}^{(2)} = \left\{
\begin{array}{ll} s_j^{(2)} e^{i
\theta_j}, \quad (n,m) \in S^{(0)} \\ 0, \quad (n,m) \in S^{(1)} \\
a \sum_l^{(2)} e^{i \theta_l}, \quad (n,m) \in S^{(2)} \\
0, \quad (n,m) \notin S^{(0)} \cup S^{(1)} \cup S^{(2)} \end{array}
\right. \qquad \psi_{n,m}^{(2)} = \left\{ \begin{array}{ll}
r_j^{(2)} e^{i \nu_j}, \quad (n,m) \in S^{(0)} \\ 0, \quad (n,m) \in S^{(1)} \\
\omega^{-2} b  \sum_l^{(2)} e^{i \nu_l}, \quad (n,m) \in S^{(2)} \\
0, \quad (n,m) \notin S^{(0)} \cup S^{(1)} \cup S^{(2)} \end{array}
\right.
\end{equation}
where $\sum_l^{(2)}$ is defined in (\ref{second-set}). The real
parameters $(s_j^{(2)},r_j^{(2)})$ satisfy an inhomogeneous system
\begin{eqnarray}
\label{inhom-eq-1} -2 (a s_j^{(2)} + \beta b r_j^{(2)}) & = & 4 +
2 \cos(\theta_{j+1} - \theta_j) + 2
\cos(\theta_{j-1} - \theta_j) + \cos(\theta_{j+2} - \theta_j), \\
\label{inhom-eq-2} -2 (\beta a s_j^{(2)} + b r_j^{(2)}) & = & 4 +
2 \cos(\nu_{j+1} - \nu_j) + 2 \cos(\nu_{j-1} - \nu_j) +
\cos(\nu_{j+2} - \nu_j).
\end{eqnarray}
When $\beta \neq 1$, the inhomogeneous system
(\ref{inhom-eq-1})--(\ref{inhom-eq-2}) has a unique solution.
Second-order corrections to the bifurcation equations are
uncoupled and have the form:
\begin{eqnarray*}
g^{(2)}_j & = & 2 \sin(\theta_j - \theta_{j+1}) + 2
\sin(\theta_j - \theta_{j-1}) + \sin(\theta_j - \theta_{j+2}), \\
h^{(2)}_j & = & 2 \sin(\nu_j - \nu_{j+1}) + 2 \sin(\nu_j -
\nu_{j-1}) + \sin(\nu_j - \nu_{j+2}),
\end{eqnarray*}
where a suitable normalization of $g_j^{(2)}$ and $h_j^{(2)}$ is
made. As a result, the Jacobian matrix computed from derivatives of
$(g_j^{(2)},h_j^{(2)})^T$ in  $(\theta_i,\nu_i)$ is block-diagonal
as ${\rm diag}({\cal M}_2,{\cal M}_2)$, where ${\cal M}_2$ is
defined in Appendix A. By Lemma 4.1 in \cite{paperII}, non-zero
eigenvalues of ${\rm diag}({\cal M}_2,{\cal M}_2)$ determine small
eigenvalues of the linearized Jacobian matrix ${\cal H}(\epsilon)$,
$$
\gamma_{1,2,3,4} = -2 \epsilon^2 + {\rm O}(\epsilon^4).
$$
Two zero eigenvalues of ${\rm diag}({\cal M}_2,{\cal M}_2)$ split
into two non-zero eigenvalues in the fourth-order LS reductions,
while two other zero eigenvalues of ${\rm diag}({\cal M}_2,{\cal
M}_2)$ persist beyond all orders due to the gauge invariance of each
component in the the coupled DNLS equations
(\ref{NLS1})--(\ref{NLS2}). Indeed, the kernel of ${\cal
H}(\epsilon)$ for $\epsilon \neq 0$ includes at least two
eigenfunctions:
\begin{equation}
\label{geometric-kernel-vector} \mbox{\boldmath $\varphi$}_{n,m} =
\left\{ \left(
\begin{array}{cc} \phi_{n,m} \\ - \bar{\phi}_{n,m} \\ 0 \\ 0 \end{array}
\right), \left(
\begin{array}{cc} 0 \\ 0 \\ \psi_{n,m} \\ - \bar{\psi}_{n,m} \end{array}
\right) \right\}, \qquad (n,m) \in \mathbb{Z}^2.
\end{equation}
In order to compute the small non-zero eigenvalues of the linearized
Jacobian matrix ${\cal H}(\epsilon)$, we use the symbolic
computation package based on Wolfram's Mathematica\footnote{The
software programs are available online at
http://dmpeli.math.mcmaster.ca/Software/LSreductions.html. For more
information on the symbolic mathematics package in which these
programs were implemented, see http://www.wolfram.com}. The
projection to the eigenspace of ${\rm diag}({\cal M}_2,{\cal M}_2)$
spanned by eigenvectors $({\bf p}_2,{\bf 0}_4)^T$ and $({\bf 0}_4,
{\bf p}_2)^T$, where ${\bf p}_2 = (-1,1,-1,1)$ and ${\bf 0}_4 =
(0,0,0,0)$, leads to the reduced eigenvalue problem (for $\omega =
1$),
\begin{eqnarray*}
\frac{-8}{1 + \beta} \left( \alpha_1 \pm \beta \alpha_2 \right) & =
& \tilde{\gamma} \alpha_1 \\
\frac{-8}{1 + \beta} \left( \pm \beta \alpha_1 + \alpha_2 \right) &
= & \tilde{\gamma} \alpha_2,
\end{eqnarray*}
where $(\alpha_1,\alpha_2)$ are coordinates of the projections,
$\tilde{\gamma} = \lim_{\epsilon \to 0} \epsilon^{-4} \gamma$, and
the upper/lower signs refer to the two coupled vortices $(1,\pm 1)$.
It is clear that the eigenvalues of the reduced eigenvalue problem
are the same for either sign and they define two small eigenvalues
of the linearized Jacobian matrix ${\cal H}(\epsilon)$ (for $\omega
= 1$):
$$
\gamma_{5} = -8 \epsilon^4 + {\rm O}(\epsilon^6), \qquad \gamma_6 =
-\frac{8(1-\beta)}{(1+\beta)} \epsilon^4 + {\rm O}(\epsilon^6).
$$

\subsection{Bifurcations of zero eigenvalues of the linearized stability problem}

We consider the eigenvalue problem (\ref{stability-vector}) in the
limit of small $\epsilon$. Let ${\cal H}(\epsilon) = {\cal H}^{(0)}
+ \epsilon {\cal H}^{(1)} + \epsilon^2 {\cal H}^{(2)} + {\rm
O}(\epsilon^3)$. The set of eigenvectors of ${\cal H}^{(0)} {\bf
e}_j = {\bf 0}$ and ${\cal H}^{(0)} {\bf f}_j = {\bf 0}$ takes the
form:
\begin{equation}
{\bf e}_j = i \left( \begin{array}{cc} e^{i \theta_j} \\ - e^{-i
\theta_j} \\ 0 \\ 0 \end{array} \right), \qquad {\bf f}_j = i \left(
\begin{array}{cc} 0 \\ 0 \\ e^{i \nu_j} \\ -e^{-i \nu_j} \end{array}
\right).
\end{equation}
The corresponding set of generalized eigenvectors of ${\cal H}^{(0)}
\hat{\bf e}_j = 2 i \sigma {\bf e}_j$ and ${\cal H}^{(0)} \hat{\bf
f}_j = 2 i \sigma {\bf f}_j$ takes the form:
\begin{equation}
\hat{\bf e}_j = \left( \begin{array}{cc} A_+ e^{i \theta_j} \\
A_+ e^{-i \theta_j} \\ B_+ e^{i \nu_j} \\ B_+ e^{-i \nu_j}
\end{array} \right), \qquad {\bf f}_j = \left(
\begin{array}{cc} A_- e^{i \theta_j} \\
A_- e^{-i \theta_j} \\ B_- e^{i \nu_j} \\ B_- e^{-i \nu_j}
\end{array} \right),
\end{equation}
where
$$
A_+ = \frac{1}{a^2 ( 1-\beta^2)}, \quad B_+ = A_- = \frac{-\beta}{a
b (1-\beta^2)}, \quad B_- = \frac{1}{b^2(1-\beta^2)}.
$$
Bifurcations of zero eigenvalues of the linear eigenvalue problem
(\ref{stability-vector}) can be computed with the extended
perturbation series expansions (\ref{vector-expansions}) for
$\phi_{n,m}(\epsilon)$ and $\psi_{n,m}(\epsilon)$ and extended
perturbation series (\ref{pertubation-eigenvectors-eigenvalues}) for
$\mbox{\boldmath $\varphi$}$ and $\lambda$, where
$$
\mbox{\boldmath $\varphi$}^{(0)} = \sum_{j = 1}^4 c_j {\bf e}_j +
\sum_{j = 1}^4 d_j {\bf f}_j, \qquad  \mbox{\boldmath
$\varphi$}^{(1)} = \frac{\lambda_1}{2} \sum_{j = 1}^4 c_j \hat{\bf
e}_j  + \frac{\lambda_1}{2} \sum_{j = 1}^4 d_j \hat{\bf f}_j +
\mbox{\boldmath $\varphi$}^{(1)}_{\rm inhom},
$$
and $\mbox{\boldmath $\varphi$}^{(1)}_{\rm inhom} = - {\cal H}^{(0)
-1} {\cal H}^{(1)} \mbox{\boldmath $\varphi$}^{(0)} = - {\cal
H}^{(1)} \mbox{\boldmath $\varphi$}^{(0)}$ is uniquely defined on
the set $S^{(1)}$. At the second-order perturbation theory, we have
the same problem (\ref{second-order-pertubation}), from which we
derive the reduced eigenvalue problem:
\begin{eqnarray}
\label{reduced-eigenvalue-vector1} {\cal M}_2 {\bf c} & = &
\frac{1}{2} \lambda_1^2 \left( A_+ {\bf c} + A_- {\bf d} \right) \\
\label{reduced-eigenvalue-vector2} {\cal M}_2 {\bf d} & = &
\frac{1}{2} \lambda_1^2 \left( B_+ {\bf c} + B_- {\bf d} \right),
\end{eqnarray}
where ${\bf c} = (c_1,c_2,c_3,c_4)^T$, ${\bf d} =
(d_1,d_2,d_3,d_4)^T$, and ${\cal M}_2$ is the same as in the scalar
case. Let $\gamma_1 = \frac{1}{2} \lambda_1^2$. The reduced
eigenvalue problem
(\ref{reduced-eigenvalue-vector1})--(\ref{reduced-eigenvalue-vector2})
has four zero roots for $\gamma_1$ and two double-degenerate
non-zero roots for $\gamma_1$, given from the quadratic equation:
\begin{equation}
(\gamma_1 + 2 a^2) (\gamma_1 + 2 b^2) = 4 a^2 b^2 \beta^2.
\end{equation}
If $\omega = 1$, such that $a^2 = b^2 = \frac{1}{1 + \beta}$, then
the two non-zero roots for $\gamma_1$ are found explicitly,
$$
\gamma_{\pm} = - \frac{2( 1 \mp |\beta|)}{1 + \beta}.
$$
By using the relation $\lambda_1 = \pm \sqrt{2 \gamma_1}$, we have
just proved that the linear eigenvalue problem
(\ref{stability-vector}) in the case $\omega = 1$ and $0 < \beta <
1$ has four small pairs of purely imaginary eigenvalues with
asymptotic approximations:
$$
\lambda_{1,2}, \lambda_{3,4} = \pm 2 i \epsilon + {\rm
O}(\epsilon^2), \qquad \lambda_{5,6},\lambda_{7,8} = \pm 2 i
\epsilon \sqrt{\frac{1-\beta}{1+\beta}} + {\rm O}(\epsilon^2).
$$
Two pairs of eigenvalues $\lambda_{5,6}$ and $\lambda_{7,8}$ become
pairs of real eigenvalues in the case $\beta > 1$. Two pairs of zero
eigenvalues of the reduced eigenvalue problem
(\ref{reduced-eigenvalue-vector1})--(\ref{reduced-eigenvalue-vector2})
split at the fourth-order LS reductions as pairs of non-zero
eigenvalues $\lambda_{9,10}$ and $\lambda_{11,12}$. Two other pairs
of zero eigenvalues persist beyond all orders for $\epsilon \neq 0$
since the geometric kernel includes two explicit solutions
(\ref{geometric-kernel-vector})) and there exists a two-parameter
solution of the inhomogeneous equation ${\cal H}(\epsilon)
\tilde{\mbox{\boldmath $\varphi$}} = 2 i \sigma \mbox{\boldmath
$\varphi$}$, where $\mbox{\boldmath $\varphi$}$ is given by
(\ref{geometric-kernel-vector}). In order to find the small non-zero
pairs of eigenvalues, we apply again the symbolic computation
package based on Wolfram's Mathematica.  The projection to the
eigenspace of ${\rm diag}({\cal M}_2,{\cal M}_2)$ spanned by
eigenvectors $({\bf p}_2,{\bf 0}_4)^T$ and $({\bf 0}_4, {\bf
p}_2)^T$ for $\lambda_1 = 0$ leads to the reduced eigenvalue problem
(for $\omega = 1$),
\begin{eqnarray*}
\frac{-8}{1 + \beta} \left( \alpha_1 \pm \beta \alpha_2 \right) & =
& \frac{1}{2 (1-\beta)} \lambda_2^2 (\alpha_1 - \beta \alpha_2)\\
\frac{-8}{1 + \beta} \left( \pm \beta \alpha_1 + \alpha_2 \right) &
= & \frac{1}{2(1-\beta)} \lambda_2^2 (-\beta \alpha_1 + \alpha_2),
\end{eqnarray*}
where $(\alpha_1,\alpha_2)$ are coordinates of the projections and
the upper/lower signs refer to the two coupled vortices $(1,\pm 1)$.
The eigenvalues of the reduced eigenvalue problem differs between
the double-charge vortex $(1,1)$ and the hidden-charge vortex
$(1,-1)$. For the double-charge vortex, the two pairs of small
eigenvalues of the linearized stability problem are purely imaginary
for any $\beta$:
$$
(1,1) : \quad \lambda_{9,10} = \pm 4 i \epsilon^2 + {\rm
O}(\epsilon^3), \qquad \lambda_{11,12} = \pm 4 i \left|
\frac{1-\beta}{1+\beta} \right| \epsilon^2 + {\rm O}(\epsilon^3).
$$
For the hidden-charge vortex, the two pairs of small eigenvalues of
the linearized stability problem are purely imaginary for $0 < \beta
< 1$ and real for $\beta > 1$:
$$
(1,-1) : \quad \lambda_{9,10},\lambda_{11,12} = \pm 4 i
\sqrt{\frac{1-\beta}{1+\beta}} \epsilon^2 + {\rm O}(\epsilon^3).
$$
We can specify precisely how many purely imaginary eigenvalues of
the linearized stability problem (\ref{stability-vector}) have
negative Krein signature. When $0 < \beta < 1$, there are eight
negative eigenvalues of ${\cal H}^{(0)}$ for the limiting solution
(\ref{zero-order-solution}) and six more small negative eigenvalues
occur for $\epsilon \neq 0$. The total number of negative
eigenvalues is reduced by two gauge symmetry constraints, such that
twelve negative eigenvalues in a constrained subspace match six
pairs of imaginary eigenvalues with negative Krein signature. When
$\beta > 1$, there are four negative eigenvalues of ${\cal H}^{(0)}$
for the limiting solution (\ref{zero-order-solution}) and five more
small negative eigenvalues occur for $\epsilon \neq 0$. The total
number of negative eigenvalues is reduced by one\footnote{When
$\beta$ is increased from $\beta < 1$ to $\beta > 1$, the Hessian
matrix related to two gauge symmetry constraints loses one positive
eigenvalue that passes through zero at $\beta = 1$ to the negative
eigenvalue for $\beta > 1$ \cite{PelKiv}.}, such that eight negative
eigenvalues in a constrained subspace match two real eigenvalues and
three pairs of imaginary eigenvalues with negative Krein signature
for the double-charge vortex and four real eigenvalues and two pairs
of imaginary eigenvalues with negative Krein signature for the
hidden-charge vortex. Therefore, the last pair of purely imaginary
eigenvalues $\lambda_{11,12}$ for the double-charge vortex has
positive Krein signature for $\beta > 1$.

\begin{figure}[tbp]
\begin{center}
\epsfxsize=7.0cm \epsffile{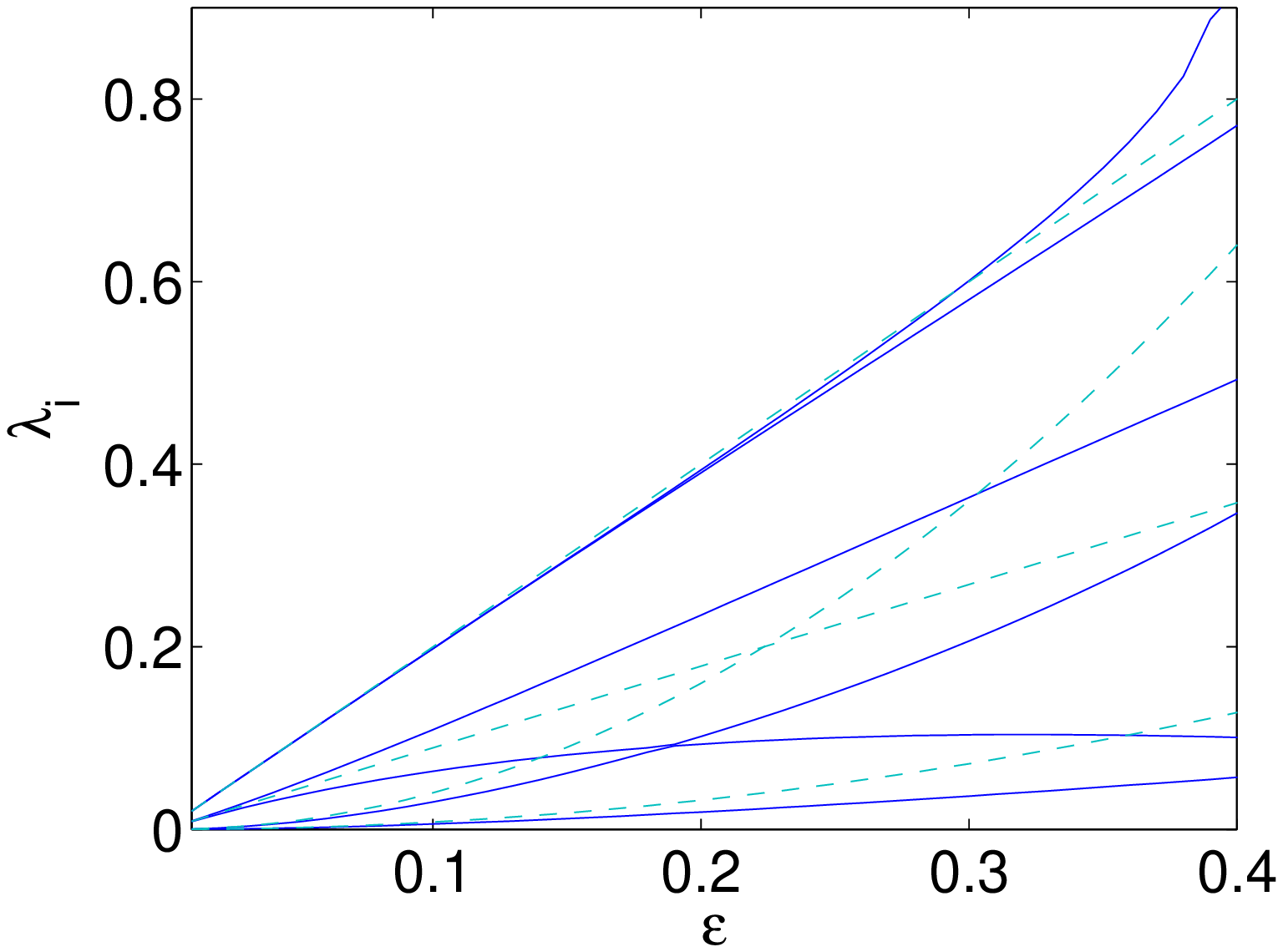} \epsfxsize=7.0cm
\epsffile{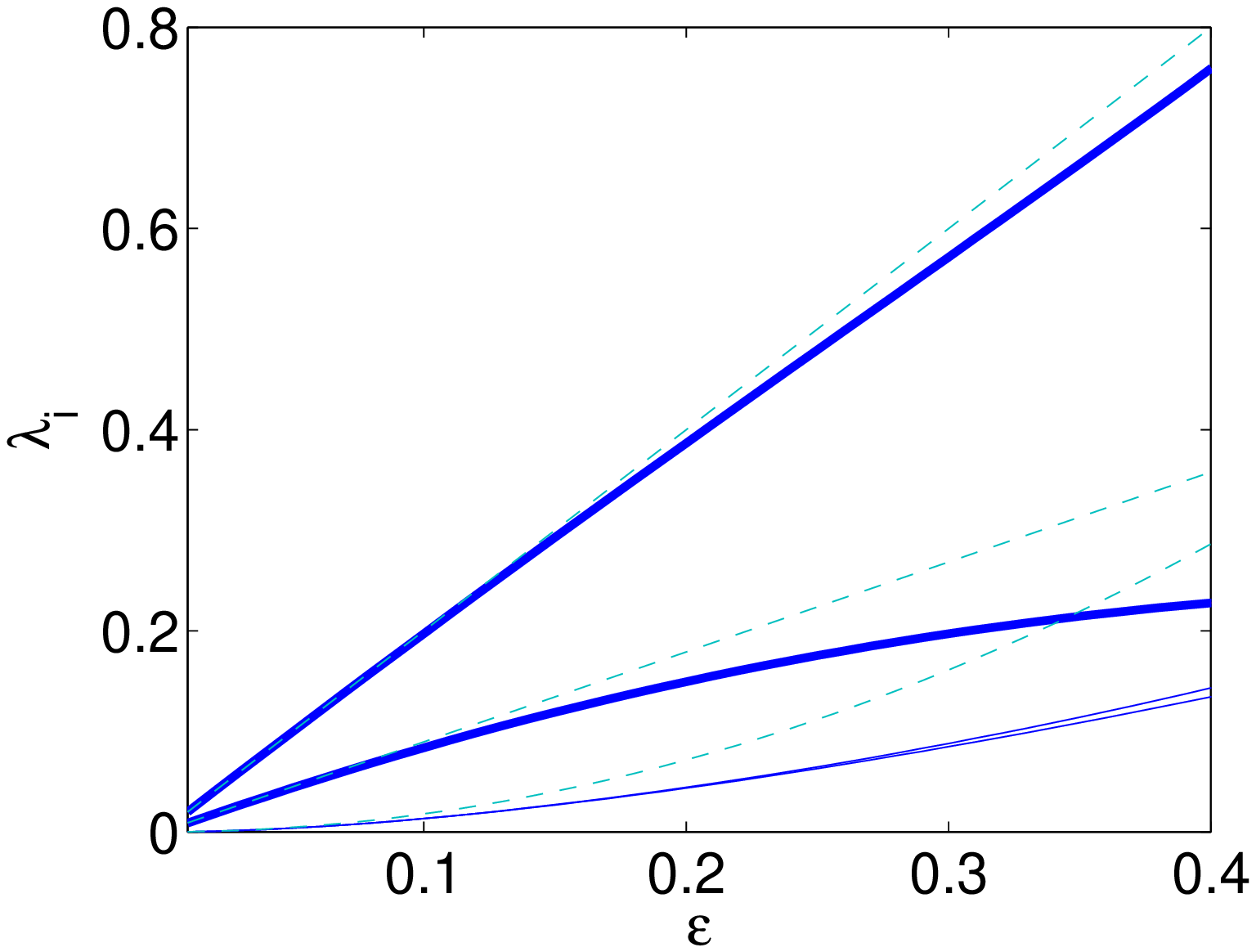} \caption{Eigenvalues of the vector vortex
cross with $\omega = 1$ and $\beta = \frac{2}{3}$ versus $\epsilon$.
Left: $(1,1)$. Right: $(1,-1)$. The solid lines show the numerical
results, while the dashed lines show the asymptotic approximations.
Bolded curves correspond to double eigenvalues (that remain
indistinguishable within the parametric window examined herein). A
good agreement is observed for $\epsilon<0.1$.} \label{fig2}
\end{center}
\end{figure}

\begin{figure}[tbp]
\begin{center}
\epsfxsize=7.0cm \epsffile{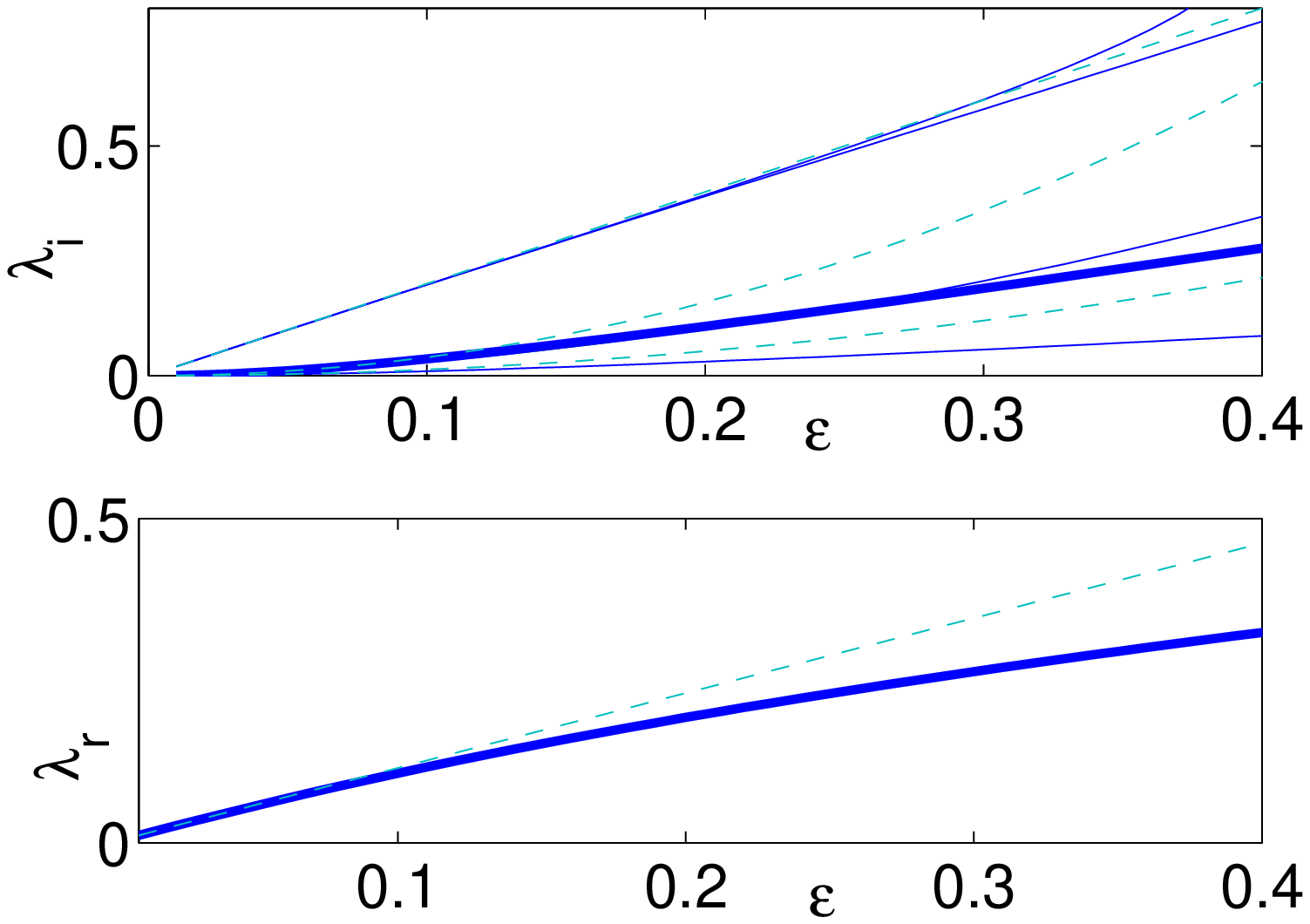} \epsfxsize=7.0cm
\epsffile{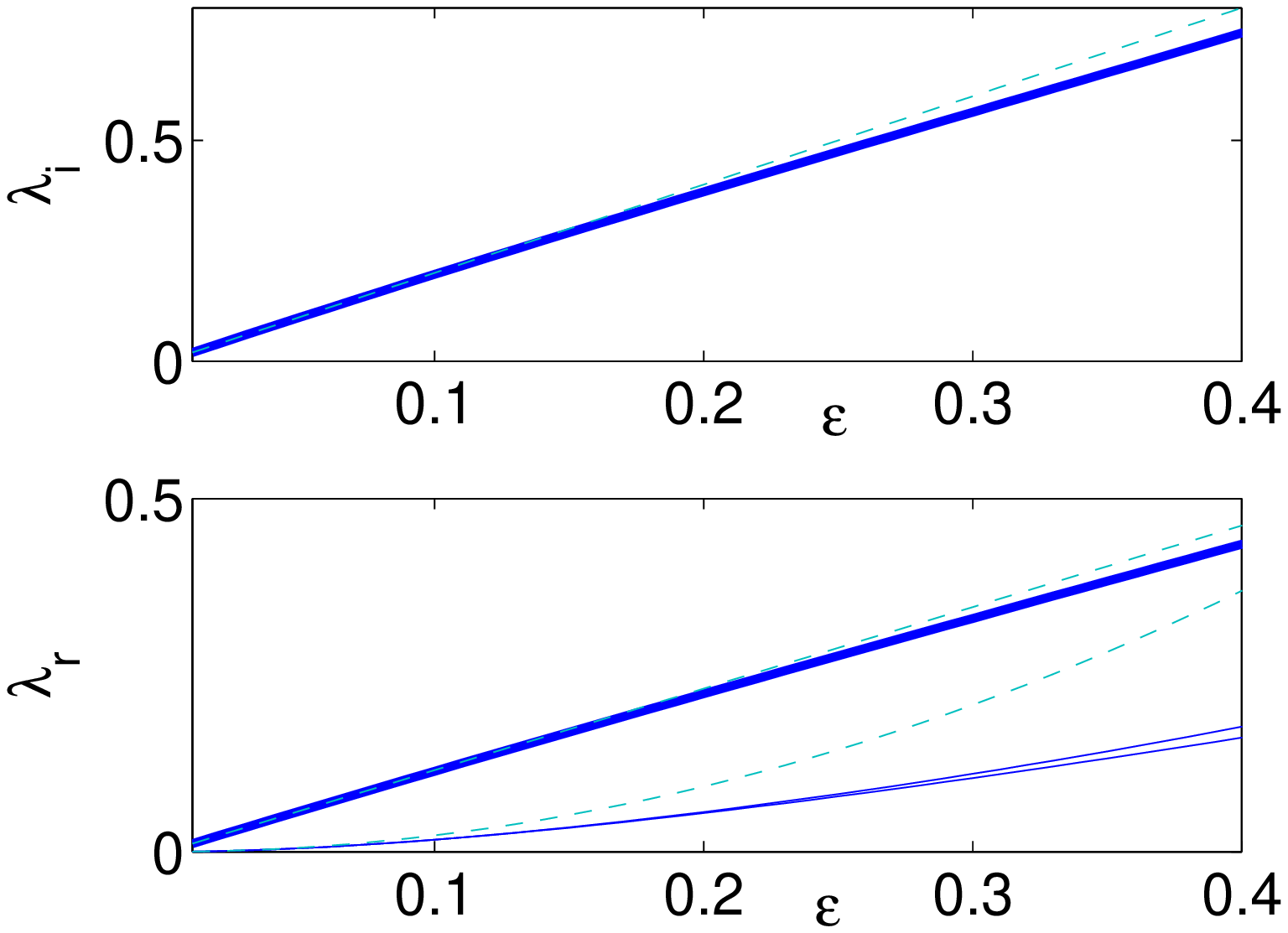} \caption{Eigenvalues of the vector vortex
cross with $\omega = 1$ and $\beta = 2$ versus $\epsilon$. Left:
$(1,1)$. Right: $(1,-1)$. The solid lines show the numerical
results, while the dashed lines show the asymptotic approximations.
Bolded curves on the left panel correspond to the real and imaginary
parts of complex eigenvalues, while bolded curves on the right panel
correspond to double eigenvalues.} \label{fig3}
\end{center}
\end{figure}

We obtain numerically small eigenvalues $\lambda$ for small values
of $\epsilon$ and $\omega = 1$. The results are shown on Figure
\ref{fig2} for $\beta = \frac{2}{3}$ and on Figure \ref{fig3} for
$\beta = 2$. The left plot corresponds to the vortex pair $(1,1)$,
while the right plot corresponds to the vortex pair $(1,-1)$. We
note that the degeneracy of the pairs $\lambda_{1,2} =
\lambda_{3,4}$ and $\lambda_{5,6} = \lambda_{7,8}$ is preserved for
the case $(1,-1)$, such that each bolded curve is double. The
degeneracy of these eigenvalues is broken for the case $(1,1)$ and
it is also broken for the pair $\lambda_{9,10} \neq \lambda_{11,12}$
for the case $(1,-1)$.

In the case of $\beta = \frac{2}{3}$, shown in Fig. \ref{fig2}, all
six pairs of neutrally stable eigenvalues bifurcate to the complex
plane for larger values of $\epsilon$ due to the Hamiltonian--Hopf
(HH) bifurcation. The first HH bifurcation happens earlier for the
case $(1,1)$ at $\epsilon \approx 0.395$, due to the broken
degeneracy between the two pairs of eigenvalues $\lambda_{1,2}$ and
$\lambda_{3,4}$. For the case $(1,-1)$, the first HH bifurcation
occurs at $\epsilon \approx 0.495$, i.e. the the hidden-charge
vortex has a larger stability window for $0 < \beta < 1$ (a similar
observation is reported for continuous systems in
\cite{desyat,fangwei}).

In the case of $\beta=2$, shown in Fig. \ref{fig3}, both cases
$(1,1)$ and $(1,-1)$ are always unstable due to the pairs of
eigenvalues $\lambda_{5,6}$ and $\lambda_{7,8}$. There are also
additional observations. In the case $(1,1)$, the pairs of double
real eigenvalues in the second-order LS reductions $\lambda_{5,6}$
and $\lambda_{7,8}$ split as a quartet of complex eigenvalues,
similarly to our computations in \cite{paperII}. Real and imaginary
parts of the quartet of complex eigenvalues are shown on Fig. 4
(left panel) by bolded curves. Only three HH bifurcations out of
four pairs of purely imaginary eigenvalues occur for larger values
of $\epsilon$. In the case $(1,-1)$, two more pairs of real
eigenvalues occur such that the hidden-charge vortex is more
unstable compared to the double-charge vortex for $\beta
> 1$. Only two HH bifurcations occur for large values of $\epsilon$.

\section{Vector vortex cross for $\beta = 1$}

In the case $\beta  = 1$, the existence domain of the coupled
vortex configurations shrinks to the line $\omega = 1$. The
zero-order solution in the anti-continuum limit is given by
(\ref{zero-order-solution}), where parameters $(a,b)$ are given by
(\ref{amplitude-parameters-symmetry}). The second-order solution
of the linear inhomogeneous system
(\ref{inhom-eq-1})--(\ref{inhom-eq-2}) with a singular matrix
exists provided that the values of $\theta_j$ and $\nu_j$ are
defined by (\ref{symmetric-vortex}). The arbitrary parameter in
the second-order solution $s_j^{(2)}$ and $r_j^{(2)}$ renormalizes
the arbitrary parameter $\delta$ in the representation
(\ref{amplitude-parameters-symmetry}).

When $\beta = \omega = 1$, the existence problem
(\ref{eq1})--(\ref{eq2}) is symmetric with respect to components
$(\phi_{n,m},\psi_{n,m})$ such that the system
(\ref{eq1})--(\ref{eq2}) can be reduced to the scalar difference
equation (\ref{scalar-eq}) with the two independent
transformations:
\begin{eqnarray*}
(1,1) & : & \quad  \phi_{n,m} = \cos \delta \; \Phi_{n,m},
\qquad \psi_{n,m} = \sin \delta \; \Phi_{n,m},\\
(1,-1) & : & \quad  \phi_{n,m} = \cos \delta \; \Phi_{n,m}, \qquad
\psi_{n,m} = \sin \delta  \; \bar{\Phi}_{n,m}.
\end{eqnarray*}
The existence result for the scalar vortex cross is formulated in
Proposition \ref{proposition-existence}. We will need the
following non-degeneracy condition for the scalar vortex cross:
\begin{equation}
\left( \sum_{(n,m) \in \mathbb{Z}^2} |\Phi_{n,m}|^2 \right)^2 \neq
\left( \sum_{(n,m) \in \mathbb{Z}^2} \Phi_{n,m}^2 \right) \left(
\sum_{(n,m) \in \mathbb{Z}^2} \bar{\Phi}_{n,m}^2 \right).
\label{non-degeneracy-condition}
\end{equation}
It is clear from the limiting solution (\ref{zero-order-solution})
that the constraint (\ref{non-degeneracy-condition}) is satisfied
for small $\epsilon$. The stability problem
(\ref{stability-vector}) is different between the cases $(1,1)$
and $(1,-1)$.

\subsection{Eigenvalues of the $(1,1)$ vortex cross}

In this case, the stability problem (\ref{stability-vector}) is
block-diagonalized under the following transformation of the four
components of the vector $\mbox{\boldmath $\varphi$}$ on the
lattice node $(n,m) \in \mathbb{Z}^2$:
\begin{eqnarray*}
(1,1) & : & \quad
\left( \begin{array}{cc} a_{n,m} \\ b_{n,m} \\ c^+_{n,m} \\
c^-_{n,m}
\end{array} \right) = \left( \begin{array}{ccccc} \cos\delta &
0 & \sin\delta & 0 \\ 0 & \cos\delta & 0 & \sin\delta
\\ -\sin\delta & 0 & \cos\delta & 0 \\ 0 & -\sin\delta & 0 &
\cos\delta \end{array} \right) \; \left( \begin{array}{cc} \varphi_1
\\ \varphi_2 \\ \varphi_3 \\ \varphi_4
\end{array} \right)_{n,m}.
\end{eqnarray*}
The components $(a_{n,m},b_{n,m})$ satisfy the linear eigenvalue
problem for scalar vortices
(\ref{stability-1})--(\ref{stability-2}). The components
$(c^+_{n,m},c^-_{n,m})$ satisfy two uncoupled self-adjoint
eigenvalue problems:
\begin{equation}
\label{uncoupled-stability} (1 - |\Phi_{n,m}|^2) c_{n,m}^{\pm} -
\epsilon \left( c_{n+1,m}^{\pm} + c_{n-1,m}^{\pm} + c_{n,m+1}^{\pm}
+ c_{n,m-1}^{\pm} \right) = \pm i \lambda c_{n,m}^{\pm}.
\end{equation}
Using the result of Proposition \ref{proposition-stability} and
equivalent computations for the uncoupled self-adjoint problems
(\ref{uncoupled-stability}), we prove the following result.

\begin{Proposition}
\label{proposition-stability-11} Let $\Phi_{n,m}(\epsilon)$, $(n,m)
\in \mathbb{Z}^2$ be a family of vortex solutions defined by
Proposition \ref{proposition-existence}. The linearized problem
(\ref{stability-vector}) in the case $\beta = 1$ for the $(1,1)$
vortex cross has zero eigenvalue of algebraic multiplicity {\em six}
and geometric multiplicity {\em five} and {\em five} small pairs of
purely imaginary eigenvalues given asymptotically by
$$
\lambda_{1,2},\lambda_{3,4} = \pm 2 i \epsilon + {\rm
O}(\epsilon^2), \quad \lambda_{5,6} = \pm 2 i \epsilon^2 + {\rm
O}(\epsilon^3), \quad \lambda_{7,8} = \pm 6 i \epsilon^2 + {\rm
O}(\epsilon^3), \quad \lambda_{9,10} = \pm 4 i \epsilon^2 + {\rm
O}(\epsilon^3).
$$
The rest of the spectrum is bounded away the origin as $\epsilon \to
0$ and it is located on the imaginary axis of $\lambda$.
\end{Proposition}

\begin{Proof}
It remains to study bifurcations of zero eigenvalues in the
self-adjoint problem (\ref{uncoupled-stability}) as $\epsilon \neq
0$.  Let us define the perturbation series for the problem
(\ref{uncoupled-stability}):
$$
{\bf c}^{\pm} = {\bf c}^{(0)} + \epsilon {\bf c}^{(1)} + \epsilon^2
{\bf c}^{(2)} + {\rm O}(\epsilon^3), \quad  \lambda = \pm i
\epsilon^2 \lambda_2 + {\rm O}(\epsilon^3).
$$
The zero-order solution is spanned by unit vectors ${\bf e}_j$ at
the $j$-th component that correspond to the node $(n,m) \in
S^{(0)}$:
$$
{\bf c}^{(0)} = \sum_{j = 1}^4 \alpha_j {\bf e}_j.
$$
The first-order correction ${\bf c}^{(1)}$ takes the form:
$$
c_{n,m}^{(1)} = \left\{
\begin{array}{ll} 0, \quad (n,m) \in S^{(0)}
\\ \sum_l^{(1)} \alpha_l, \quad (n,m) \in S^{(1)} \\ 0, \quad (n,m)
\notin S^{(0)} \cup S^{(1)} \end{array} \right.
$$
where the sum $\sum_l^{(1)}$ is defined in (\ref{first-set}). At the
second-order in $\epsilon$, we find a set of non-trivial equations
at the nodes $(n,m) \in S^{(0)}$:
$$
\alpha_j + \alpha_{j+2} + 2 (\alpha_{j+1} + \alpha_{j-1}) =
\lambda_2 \alpha_j, \qquad j = 1,2,3,4.
$$
The reduced eigenvalue problem has a double zero eigenvalue and two
non-zero eigenvalues $-2$ and $6$. Two zero eigenvalues of the
problem (\ref{uncoupled-stability}) persist at all orders of
$\epsilon$, because of the exact solutions: $c_{n,m}^{\pm} =
\Phi_{n,m}$ and $c_{n,m}^{\pm} = \bar{\Phi}_{n,m}$.
\end{Proof}

We note that the pairs of eigenvalues $\lambda_{1,2}$,
$\lambda_{3,4}$, and $\lambda_{9,10}$ continue the eigenvalues of
the vortex cross $(1,1)$ from $\beta \neq 1$ to $\beta = 1$. The
pairs of eigenvalues $\lambda_{5,6}$ and $\lambda_{7,8}$ match
with the zero ${\rm O}(\epsilon)$ corrections to the corresponding
eigenvalues of the vortex cross $(1,1)$ for $\beta \neq 1$.
Finally, the pair of non-zero eigenvalues $\lambda_{11,12}$  for
$\beta \neq 1$ is forced to remain at the origin for $\beta = 1$
due to the polarization-rotation symmetry.

We can now specify how many purely imaginary eigenvalues $\lambda$
have negative Krein signature. When $\beta = 1$, there are four
negative and twelve zero eigenvalues of ${\cal H}^{(0)}$ for the
limiting solution (\ref{zero-order-solution}). Out of the twelve
zero eigenvalues, three small negative eigenvalues bifurcate in
the subspace for components $(a_{n,m},b_{n,m})$, two small
positive and two small negative eigenvalues bifurcate in the
subspace for components $(c_{n,m}^+,c_{n,m}^-)$ and five
eigenvalues remain at zero as $\epsilon \neq 0$. The total number
of negative eigenvalues is reduced by one symmetry
constraint\footnote{The Hessian matrix related to two gauge
symmetry constraints has a zero eigenvalue for $\beta = 1$, while
only positive eigenvalues are counted in a reduction of the
negative index of ${\cal H}(\epsilon)$.}, such that eight negative
eigenvalues in a constrained subspace match four pairs of
imaginary eigenvalues with negative Krein signature. The only pair
of purely imaginary eigenvalues with positive Krein signature is
the pair $\lambda_{5,6}$ that is related to the two small positive
eigenvalues in the subspace for components
$(c_{n,m}^+,c_{n,m}^-)$.

\subsection{Eigenvalues of the $(1,-1)$ vortex cross}

Since the stability problem (\ref{stability-vector}) has no
block-diagonalization for the (1,-1) vortex cross, the results of
the second-order LS reductions give only two pairs of purely
imaginary eigenvalues $\lambda_{1,2}$ and $\lambda_{3,4}$. We
shall study the eigenvalues of the fourth-order LS reduction by
using the symbolic computation package based on Wolfram's
Mathematica. In order to prepare for symbolic computations, we
note that the eigenvalues of ${\cal H}(\epsilon)$ in the case
$(1,-1)$ are exactly the same as eigenvalues of ${\cal
H}(\epsilon)$ in the case $(1,1)$, due to the equivalent
transformation of the vector $\mbox{\boldmath $\varphi$}$ in the
eigenvalue problem ${\cal H}(\epsilon) \mbox{\boldmath $\varphi$}
= \gamma \mbox{\boldmath $\varphi$}$:
\begin{eqnarray*}
(1,-1) & : & \quad
\left( \begin{array}{cc} a_{n,m} \\ b_{n,m} \\ c^+_{n,m} \\
c^-_{n,m}
\end{array} \right) = \left( \begin{array}{ccccc} \cos\delta &
0 & 0 & \sin\delta  \\ 0 & \cos\delta & \sin\delta & 0
\\ 0 & -\sin\delta & \cos\delta & 0 \\ -\sin\delta & 0 & 0 &
\cos\delta \end{array} \right) \; \left( \begin{array}{cc} \varphi_1
\\ \varphi_2 \\ \varphi_3 \\ \varphi_4
\end{array} \right)_{n,m}.
\end{eqnarray*}
As a result of this transformation, we immediately find the {\em
five-dimensional} kernel of ${\cal H}(\epsilon)$ for $\epsilon
\neq 0$, which can be spanned as follows: {\small
\begin{equation}
\label{geometric-kernel-Manakov} \mbox{\boldmath $\varphi$}_{n,m}
= \left\{ \left(
\begin{array}{cc} \cos \delta \; \Phi_{n,m} \\ - \cos \delta \;\bar{\Phi}_{n,m} \\
-\sin \delta \; \bar{\Phi}_{n,m} \\ \sin \delta \; \Phi_{n,m}
\end{array} \right), \left(
\begin{array}{cc} -\sin \delta \; \Phi_{n,m} \\ 0 \\
0 \\ \cos \delta \; \Phi_{n,m}
\end{array} \right), \left(
\begin{array}{cc} 0 \\ - \sin \delta \;\Phi_{n,m} \\
\cos \delta \; \Phi_{n,m} \\ 0
\end{array} \right), \left(
\begin{array}{cc} -\sin \delta \; \bar{\Phi}_{n,m} \\ 0 \\
0 \\ \cos \delta \; \bar{\Phi}_{n,m}
\end{array} \right), \left(
\begin{array}{cc} 0 \\ - \sin \delta \;\bar{\Phi}_{n,m} \\
\cos \delta \; \bar{\Phi}_{n,m} \\ 0
\end{array} \right),  \right\},
\end{equation}
}for $(n,m) \in \mathbb{Z}^2$. Algebraic multiplicity of zero
eigenvalue for $\epsilon \neq 0$ is defined by the solution of the
inhomogeneous equation ${\cal H}(\epsilon) \tilde{\mbox{\boldmath
$\varphi$}} = 2 i \sigma \mbox{\boldmath $\varphi$}$, which is
equivalent to the projection equations
$$
\sum_{(n,m) \in \mathbb{Z}^2} \langle \mbox{\boldmath
$\varphi$}_j, \sigma \mbox{\boldmath $\varphi$} \rangle = 0,
\qquad j = 1,2,3,4,5,
$$
where $\mbox{\boldmath $\varphi$}$ is spanned by five eigenvectors
$\mbox{\boldmath $\varphi$}_j$ in the decomposition
(\ref{geometric-kernel-vector}). Solving this system of linear
equations, we have found under the non-degeneracy condition
(\ref{non-degeneracy-condition}) that there is a one-parameter
solution of the inhomogeneous system for $\delta \neq
\frac{\pi}{4}$ and a three-parameter solution for $\delta =
\frac{\pi}{4}$. Thus, the zero eigenvalue has algebraic
multiplicity {\em six} for $\delta \neq \frac{\pi}{4}$ and {\em
eight} for $\delta = \frac{\pi}{4}$.

In the limit $\epsilon = 0$, when ${\cal H}^{(0)} = {\cal H}(0)$,
we construct explicitly three sets of linearly independent
eigenvectors of ${\cal H}^{(0)}$:
\begin{equation}
{\bf e}_j = i \left( \begin{array}{cc} \cos \delta \; e^{i \theta_j} \\
- \cos \delta \; e^{-i \theta_j} \\ \sin \delta \; e^{-i \theta_j} \\
-\sin \delta \; e^{i \theta_j} \end{array} \right), \quad {\bf
f}_j^+ = i \left( \begin{array}{cc} \cos \delta \; e^{i \theta_j} \\
- \cos \delta \; e^{-i \theta_j} \\ -\sin \delta \; e^{-i \theta_j} \\
\sin \delta \; e^{i \theta_j}\end{array} \right), \quad {\bf
f}_j^- = \left(
\begin{array}{cc} \sin \delta \; e^{i \theta_j} \\
\sin \delta \; e^{-i \theta_j} \\ -\cos \delta \; e^{-i \theta_j} \\
-\cos \delta \; e^{i \theta_j}\end{array} \right),
\end{equation}
Only the set of eigenvectors ${\bf e}_j$ generates the set of
generalized eigenvectors of the problem ${\cal H}^{(0)} \hat{\bf
e}_j = 2 i \sigma {\bf e}_j$, where
\begin{equation}
\hat{\bf e}_j = \left( \begin{array}{cc} \cos \delta \; e^{i
\theta_j}\\ \cos \delta \; e^{-i \theta_j} \\ \sin \delta \; e^{-i \theta_j} \\
\sin \delta \; e^{i \theta_j} \end{array} \right).
\end{equation}
Thus, the zero eigenvalue of ${\cal H}^{(0)}$ has algebraic
multiplicity {\em sixteen} and geometric multiplicity {\em
twelve}. Two pairs of purely imaginary eigenvalues of negative
Krein signatures bifurcate at the second-order LS reductions as
$$
\lambda_{1,2},\lambda_{3,4} = \pm 2 i \epsilon + {\rm
O}(\epsilon^2).
$$
In order to study bifurcations of non-zero eigenvalues at the
fourth-order LS reductions, we consider the extended perturbation
series (\ref{pertubation-eigenvectors-eigenvalues}) for
$\mbox{\boldmath $\varphi$}$ and $\lambda$ with $\lambda_1 = 0$ and
\begin{equation}
\label{decomposition-zero-solution} \mbox{\boldmath $\varphi$}^{(0)}
= \sum_{j = 1}^4 c_j {\bf e}_j + \sum_{j = 1}^4 d_j^+ {\bf f}_j^+ +
\sum_{j = 1}^4 d_j^- {\bf f}_j^-.
\end{equation}
Performing computations symbolically, we have twelve homogeneous
equations at the order of ${\rm O}(\epsilon^2)$ for twelve variables
$(c_j,d_j^+,d_j^-)$, $j = 1,2,3,4$, which can be converted and
simplified to the following determinant equation:
$$
\gamma_2^2 + 4 ( 1 + 4 \cos 4 \delta ) \gamma_2 + 36 = 0,
$$
where $\gamma_2 = \frac{1}{2} \lambda_2^2$. By using the inverse
relation $\lambda_2 = \pm \sqrt{2 \gamma_2}$ and finding the roots
for $\gamma_2$ explicitly, we obtain four small pairs of
eigenvalues with asymptotic approximations:
\begin{eqnarray*}
\lambda_{5,6} & = & \pm 2 i \epsilon^2 \sqrt{1 + 4 \cos 4 \delta -
\sqrt{8(\cos 4 \delta + \cos 8 \delta)}} + {\rm O}(\epsilon^3), \\
\lambda_{7,8} & = & \pm 2 i \epsilon^2 \sqrt{1 + 4 \cos 4 \delta +
\sqrt{8(\cos 4 \delta + \cos 8 \delta)}} + {\rm O}(\epsilon^3).
\end{eqnarray*}
When $\delta = 0$ or $\delta = \frac{\pi}{2}$, we obtain the same
pairs of purely imaginary eigenvalues as in the case $(1,1)$ (see
Proposition \ref{proposition-stability-11}). When $\delta =
\frac{\pi}{4}$, we obtain two degenerate pairs of real eigenvalues
$$
\lambda_{5,6},\lambda_{7,8} = \pm 2 \sqrt{3} \epsilon^2 + {\rm
O}(\epsilon^3).
$$
The instability domain is found analytically from the condition
that complex-valued roots for $\gamma_2$ coalesce and become a
double negative root. This happens when $\cos(4 \delta) + \cos(8
\delta) = 0$, which is solved on the interval $\delta \in \left[
0,\frac{\pi}{2}\right]$ at $\delta = \frac{\pi}{12}$ and $\delta =
\frac{5 \pi}{12}$\footnote{Another solution exists at $\delta =
\frac{\pi}{4}$ but it corresponds to the case when complex-valued
roots coalesce and become a double positive root for $\gamma_2$.}
Thus, the instability domain of the $(1,-1)$ vortex cross in the
case $\beta = 1$ is bounded by the interval $\delta \in \left(
\frac{\pi}{12},\frac{5 \pi}{12} \right)$.

In order to capture the remaining pair of non-zero eigenvalues
$\lambda_{9,10}$, we shall reorder the perturbation series
expansions and to move the last two sums in the decomposition
(\ref{decomposition-zero-solution}) to the order of ${\rm
O}(\epsilon^2)$, while the coefficients of the vector ${\bf c} =
(c_1,c_2,c_3,c_4)^T$ should be projected to the vector ${\bf p}_2
= (-1,1,-1,1)$ of the kernel of ${\cal M}_2$, such that ${\bf c} =
x_1 {\bf p}_2$. Performing computations symbolically, we have
twelve homogeneous equations at the order of ${\rm O}(\epsilon^4)$
for eight variables in the vectors ${\bf d}^+$ and ${\bf d}_j^-$
and the coordinate $x_1$. The homogeneous system is satisfied with
the choice ${\bf d}^+ = {\bf 0}$ and ${\bf d}^- = x_2 {\bf p}_2$,
where $x_2$ is another coordinate. The coordinates $(x_1,x_2)$
solve a homogeneous system with the determinant equation
$\lambda_2^2 = - 16 \cos^2(2\delta)$. Therefore, a small pair of
purely imaginary eigenvalues of negative Krein signatures has the
asymptotic approximation:
\begin{eqnarray*}
\lambda_{9,10} = \pm 4 i \epsilon^2 \cos(2\delta) + {\rm
O}(\epsilon^3).
\end{eqnarray*}
When $\delta = 0$ and $\delta = \frac{\pi}{2}$, the pair
$\lambda_{9,10}$ matches to that in the case $(1,1)$ (see
Proposition \ref{proposition-stability-11}). When $\delta =
\frac{\pi}{4}$, the pair remains at the origin as it follows from
the study of algebraic multiplicity of zero eigenvalue. According
to the count of negative eigenvalues, the total number of negative
eigenvalues of ${\cal H}(\epsilon)$ for small $\epsilon$ reduced
by one symmetry constraint is eight. These eigenvalues match two
pairs of imaginary eigenvalues $\lambda_{1,2}$ and $\lambda_{3,4}$
and two real positive eigenvalues $\lambda_{5,6}$ and
$\lambda_{7,8}$\footnote{Eigenvalues $\lambda_{5,6}$ and
$\lambda_{7,8}$ are real only in the case $\delta =
\frac{\pi}{4}$. For other values of $\delta$, these eigenvalues
are either complex-valued or purely imaginary. The count is not
affected, since two real eigenvalues are equivalent to four
complex eigenvalues which may coalesce due to the inverse
Hamilton--Hopf bifurcation to two pairs of purely imaginary
eigenvalues with positive and negative Krein signatures.}.

\begin{figure}[tbp]
\begin{center}
\epsfxsize=7.0cm \epsffile{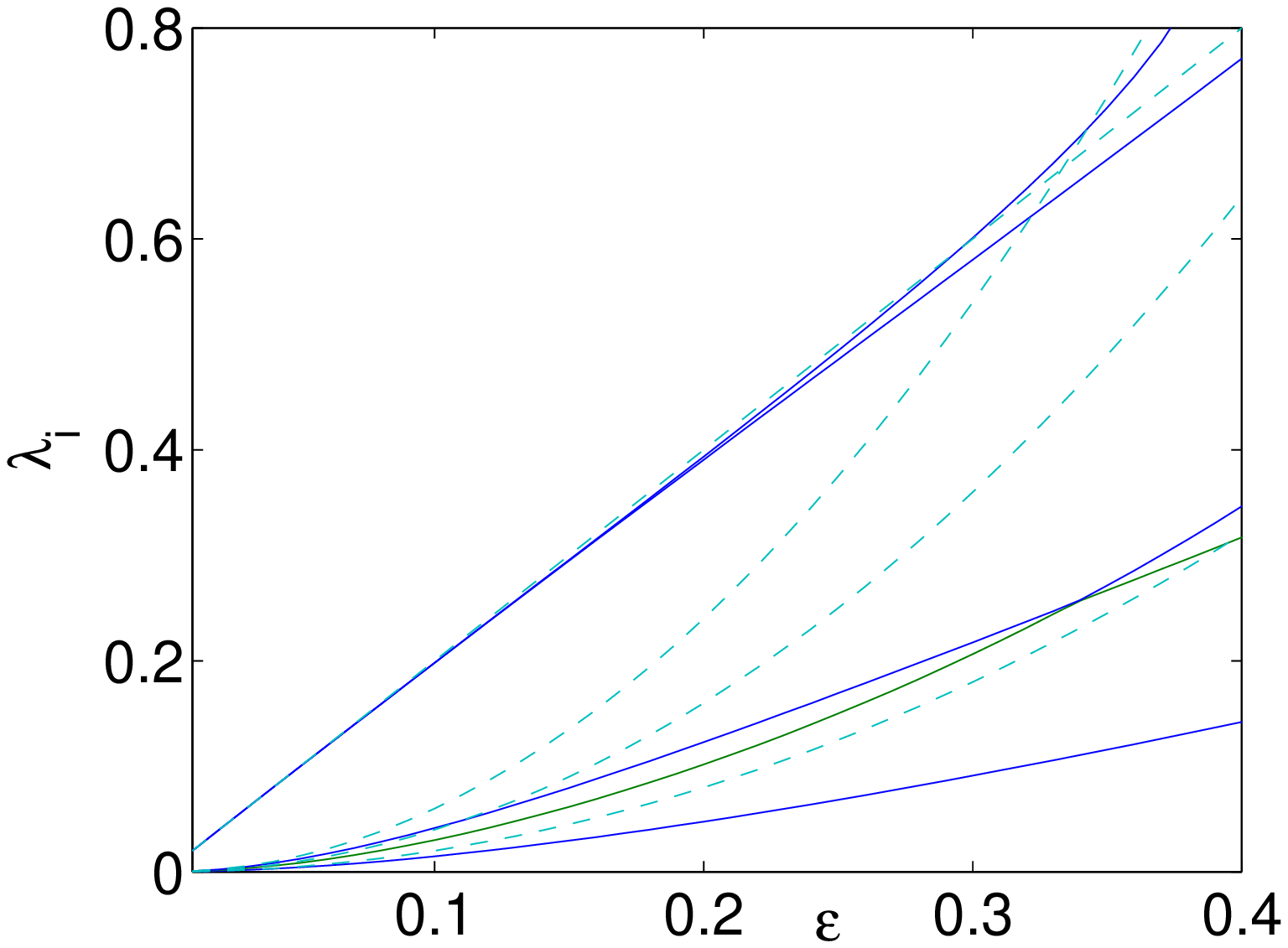} \epsfxsize=7.0cm
\epsffile{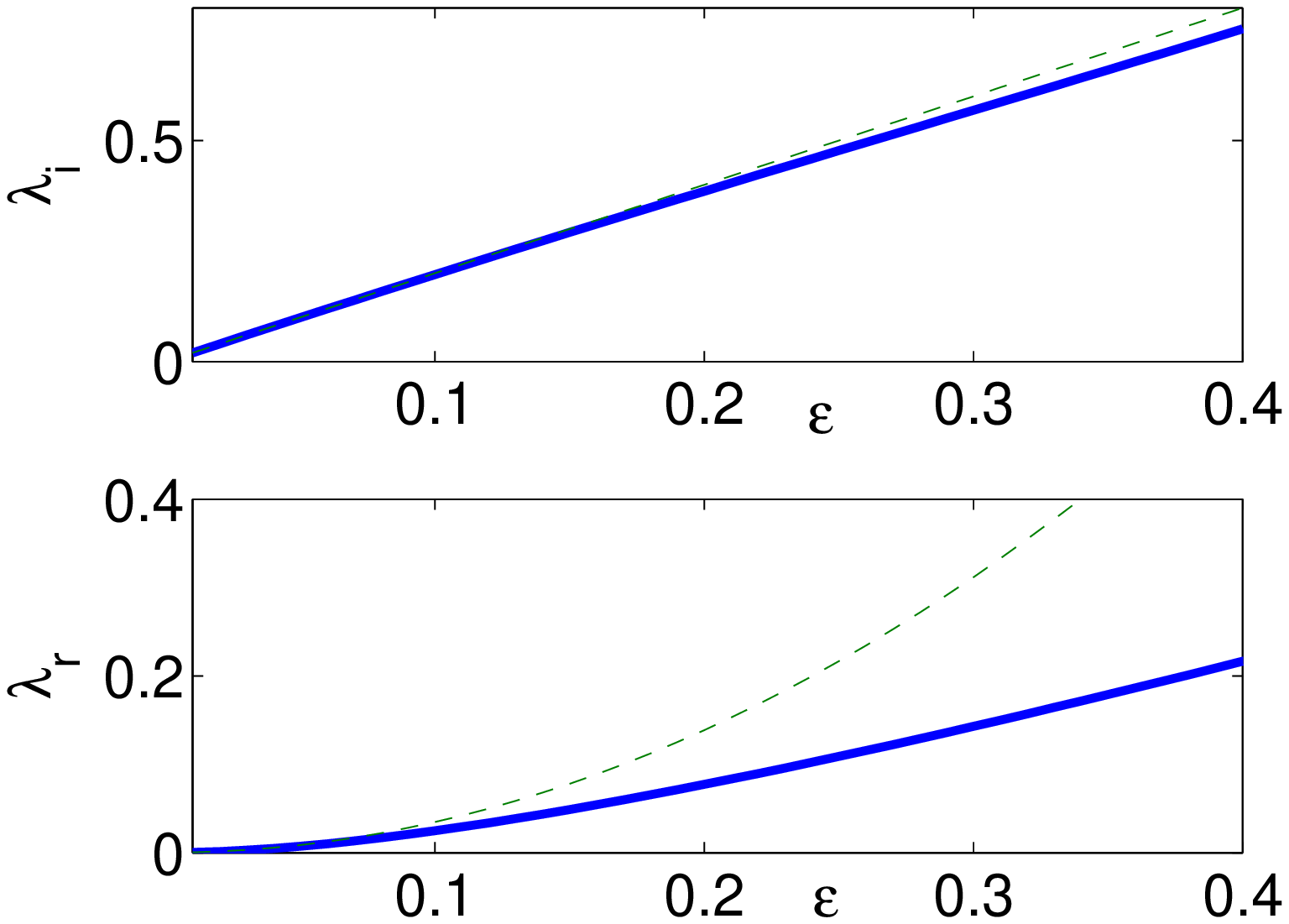} \caption{Eigenvalues of the vector
vortex cross with $\omega = \beta = 1$ and $\delta =
\frac{\pi}{4}$ versus $\epsilon$. Left: $(1,1)$. Right: $(1,-1)$.
The solid lines show the full numerical results, while the dashed
lines show the asymptotic approximations. Bolded curves show
double eigenvalues.} \label{fig4}
\end{center}
\end{figure}

Asymptotic and numerical approximations of small eigenvalues
$\lambda$ for small values of $\epsilon$ for $\omega = \beta = 1$
and $\delta = \frac{\pi}{4}$ are shown on Figure \ref{fig4}. The
left plot corresponds to the vortex pair $(1,1)$, while the right
plot corresponds to the vortex pair $(1,-1)$. We can see that the
$(1,1)$ vortex cross is linearly stable in the anti-continuum
limit, according to the results of Proposition
\ref{proposition-stability-11}. On the other hand, the $(1,-1)$
vortex cross become unstable because of the a double pairs of real
eigenvalues $\lambda_{5,6} = \lambda_{7,8}$. The other double pair
of purely imaginary eigenvalues remains double for all $\epsilon >
0$, such that $\lambda_{1,2} = \lambda_{3,4}$. Therefore, the
stability changes drastically in the case of the discrete Manakov
system (that is the coupled DNLS system for $\beta = 1$): the
$(1,1)$ vortex cross is stable near the anti-continuum limit while
the $(1,-1)$ vortex cross is linearly unstable.

\section{Conclusions}

We have examined analytically and numerically the existence and
stability of vortex cross configurations in the single-component and
two-component DNLS equations. We have used the Lyapunov-Schmidt
theory, to obtain the bifurcation functions and the solvability
conditions that allow persistence of such configurations near the
anti-continuum limit. Additionally, the theory gives analytical
expressions for eigenvalues of the linearized stability problem as
functions of the system parameters (namely, the coupling between
adjacent lattice nodes $\epsilon \geq 0$ and the coupling between
the two components $\beta \geq 0$).

One of the interesting recent experimental developments in the
setting of BECs concerns the experimental and theoretical studies of
spin-1 (or 3-component) states \cite{higbie,higbie2}. This feature,
along with similar possibilities that could be realized in optical
settings, render desirable a general theory for interactions of
multiple components. Such studies are currently in progress.

\appendix

\section{Continuation of the single-component vortex cross}

We apply the algorithm of Lyapunov--Schmidt (LS) reductions (see
\cite{paperII} for details) and compute the first few terms of the
perturbation series expansions:
\begin{equation}
\label{appendix-expansions} \Phi_{n,m}(\epsilon) = \sum_{k =
0}^{\infty} \epsilon^k \Phi_{n,m}^{(k)}.
\end{equation}
The zero-order solution $\Phi_{n,m}^{(0)}$ is given by
(\ref{scalar-vortex-0}). The first-order correction is obtained in
the explicit form:
\begin{equation}
\Phi_{n,m}^{(1)} = \left\{ \begin{array}{ll} 0, \quad (n,m) \in
S^{(0)}
\\ \sum_l^{(1)} e^{i \theta_l}, \quad (n,m) \in S^{(1)} \\ 0, \quad (n,m)
\notin S^{(0)} \cup S^{(1)} \end{array} \right.
\end{equation}
where $S^{(1)}$ is the set of adjacent nodes to the set $S^{(0)}$
and $\sum_l^{(1)} e^{i \theta_l}$ is a schematic notation for the
following solution:
\begin{equation}
\label{first-set} \sum_l^{(1)} e^{i \theta_l} = \left\{
\begin{array}{ll} e^{i
\theta_1} + e^{i \theta_2} + e^{i \theta_3} + e^{i \theta_4}, \quad (n,m) = (0,0) \\
e^{i \theta_j} + e^{\theta_{j+1}}, \quad (n,m) = \{
(-1,-1);(1,-1);(1,1);(-1,1)\} \\
e^{i \theta_j}, \quad (n,m) = \{ (-2,0);(0,-2);(2,0);(0,2)\}
\end{array} \right.
\end{equation}
The index $j$ enumerates nodes in the set $S^{(0)}$ that are
adjacent to the nodes in the set $S^{(1)}$ listed in the figured
brackets of (\ref{first-set}). No non-trivial bifurcation
equations arise at the first-order reductions, i.e. the
first-order correction to the bifurcation function ${\bf
g}^{(1)}(\mbox{\boldmath $\theta$})$ is zero, where
$\mbox{\boldmath $\theta$} =
(\theta_1,\theta_2,\theta_3,\theta_4)$ and notations of
\cite{paperII} are used. The second-order correction is found in
the form:
\begin{equation}
\Phi_{n,m}^{(2)} = \left\{ \begin{array}{ll} s_j^{(2)}  e^{i
\theta_j}, \quad (n,m) \in S^{(0)} \\ 0, \quad (n,m) \in S^{(1)} \\
\sum_l^{(2)} e^{i \theta_l}, \quad (n,m) \in S^{(2)} \\
0, \quad (n,m) \notin S^{(0)} \cup S^{(1)} \cup S^{(2)} \end{array}
\right.
\end{equation}
where
$$
-2 s_j^{(2)} = 4 + 2 \cos(\theta_{j+1} - \theta_j) + 2
\cos(\theta_{j-1} - \theta_j) + \cos(\theta_{j+2} - \theta_j).
$$
The set $S^{(2)}$ contains outward adjacent nodes to the set
$S^{(1)} \backslash \{(0,0)\}$ and $\sum_l^{(2)} e^{i \theta_l}$
is a schematic notation for the following solution:
\begin{equation}
\label{second-set} \sum_l^{(2)} e^{i \theta_l} = \left\{
\begin{array}{ll} 2 e^{i \theta_j} + e^{i \theta_{j+1}}, \quad (n,m)
= \{
(-2,-1);(1,-2);(2,1);(-1,2)\} \\
e^{i \theta_j} + 2 e^{\theta_{j+1}}, \quad (n,m) = \{
(-1,-2);(2,-1);(1,2);(-2,1)\} \\
e^{i \theta_j}, \quad (n,m) = \{ (-3,0);(0,-3);(3,0);(0,3)\}
\end{array} \right.
\end{equation}
The second-order corrections to the bifurcation function take the
form:
\begin{equation}
g^{(2)}_j = 2 \sin(\theta_j - \theta_{j+1}) + 2 \sin(\theta_j -
\theta_{j-1}) + \sin(\theta_j - \theta_{j+2}), \qquad j = 1,2,3,4.
\end{equation}
The bifurcation equations ${\bf g}^{(2)}(\mbox{\boldmath $\theta$})
= {\bf 0}$ are satisfied with the one-parameter family of asymmetric
vortices:
\begin{equation}
\label{asymmetric-vortices} \theta_1 = 0, \quad \theta_2 = \theta,
\quad \theta_3 = \pi, \quad \theta_4 = \pi + \theta,
\end{equation}
where $\theta \in (0,\pi)$. When $\theta = \frac{\pi}{2}$, the
family (\ref{asymmetric-vortices}) reduces to the vortex cross
configuration (\ref{symmetric-vortex}). The Jacobian matrix ${\cal
M}_2$ of the second-order bifurcation function ${\bf
g}^{(2)}(\mbox{\boldmath $\theta$})$ is obtained by differentiation
of ${\bf g}^{(2)}$ in $\mbox{\boldmath $\theta$}$. At the family of
asymmetric vortices (\ref{asymmetric-vortices}), the Jacobian matrix
${\cal M}_2$ takes the form:
$$
{\cal M}_2 = \left( \begin{array}{ccccc} -1 & - 2 \cos \theta & 1 &
2 \cos \theta \\ -2 \cos \theta & -1 & 2 \cos \theta & 1 \\ 1 & 2
\cos \theta & -1 & -2 \cos \theta \\ 2 \cos \theta & 1 & -2\cos
\theta & -1 \end{array} \right).
$$
It has two zero eigenvalues and two non-zero eigenvalues $-2 \pm 4
\cos \theta$. In the case of the vortex cross ($\theta =
\frac{\pi}{2}$), it has two zero eigenvalues and two negative
eigenvalues $-2$. The third-order correction satisfies the
inhomogeneous equation,
$$
(1 - 2 |\Phi_{n,m}^{(0)}|^2 ) \Phi_{n,m}^{(3)} - \Phi_{n,m}^{(0)2}
\bar{\Phi}_{n,m}^{(3)} = \Phi_{n+1,m}^{(2)} + \Phi_{n-1,m}^{(2)} +
\Phi_{n,m+1}^{(2)} + \Phi_{n,m-1}^{(2)} + |\Phi_{n,m}^{(1)}|^2
\Phi_{n,m}^{(1)}
$$
where we have shorten nonlinear terms, since $\Phi_{n,m}^{(0)}
\Phi_{n,m}^{(1)} = \Phi_{n,m}^{(1)} \Phi_{n,m}^{(2)} = 0$ for all
$(n,m) \in \mathbb{Z}^2$. The third-order correction is found in the
form:
\begin{equation}
\Phi_{n,m}^{(3)} = \left\{ \begin{array}{ll} 0, \quad (n,m) \in S^{(0)} \\
|\Phi_{n,m}^{(1)}|^2 \Phi_{n,m}^{(1)} + \sum_l^{(1)} s_l^{(2)} e^{i
\theta_l}
+ \sum_l^{(1,3)} e^{i \theta_l}, \quad (n,m) \in S^{(1)} \\
0, \quad (n,m) \in S^{(2)} \\
\sum_l^{(3)} e^{i \theta_l}, \quad (n,m) \in S^{(3)} \\
0, \quad (n,m) \notin S^{(0)} \cup S^{(1)} \cup S^{(2)} \cup S^{(3)}
\end{array} \right.
\end{equation}
where the sum $\sum_l^{(1)} s_l^{(2)} e^{i \theta_l}$ is defined
similarly to the sum (\ref{first-set}), the sum $\sum_l^{(3)} e^{i
\theta_l}$ is not used for further computations, and the sum
$\sum_l^{(1,3)} e^{i \theta_l}$ is defined as follows:
\begin{equation}
\label{first-set-modified} \sum_l^{(1,3)} e^{i \theta_l} = \left\{
\begin{array}{ll} 3 e^{i \theta_j} + 3 e^{\theta_{j+1}}, \quad (n,m) = \{
(-1,-1);(1,-1);(1,1);(-1,1)\} \\
5 e^{i \theta_j} + e^{i \theta_{j+1}} + e^{i \theta_{j-1}}, \quad
(n,m) = \{ (-2,0);(0,-2);(2,0);(0,2)\} \end{array} \right.
\end{equation}
No non-trivial bifurcation equations arise at the third-order
reductions, i.e. ${\bf g}^{(3)}(\mbox{\boldmath $\theta$}) = {\bf
0}$. The fourth-order correction satisfies the inhomogeneous
equation, {\small
\begin{eqnarray*}
(1 - 2 |\Phi_{n,m}^{(0)}|^2 ) \Phi_{n,m}^{(4)} - \Phi_{n,m}^{(0)2}
\bar{\Phi}_{n,m}^{(4)} =  2 |\Phi_{n,m}^{(2)}|^2 \Phi_{n,m}^{(0)}
+ \Phi_{n,m}^{(2)2} \bar{\Phi}_{n,m}^{(0)} + \Phi_{n+1,m}^{(3)} +
\Phi_{n-1,m}^{(3)} + \Phi_{n,m+1}^{(3)} + \Phi_{n,m-1}^{(3)}.
\end{eqnarray*}
} Solving the inhomogeneous equation for the third-order
corrections, we obtain the bifurcation equations at the fourth
order of LS reductions in the form: {\small
\begin{eqnarray*}
\nonumber g^{(4)}_j & = & \left( 4 + 2
\cos(\theta_{j+2}-\theta_{j+1}) + 2 \cos(\theta_j - \theta_{j+1}) +
\cos(\theta_{j-1} - \theta_{j+1}) \right) \sin(\theta_{j+1} -
\theta_j) \\
\nonumber & + & \left( 4 + 2 \cos(\theta_{j-2}-\theta_{j-1}) + 2
\cos(\theta_j - \theta_{j-1}) + \cos(\theta_{j+1} - \theta_{j-1})
\right) \sin(\theta_{j-1} - \theta_j) \\
\nonumber & + & \frac{1}{2} \left( 4 + 2
\cos(\theta_{j-1}-\theta_{j-2}) + 2 \cos(\theta_{j+1} -
\theta_{j+2}) + \cos(\theta_j - \theta_{j+2})
\right) \sin(\theta_{j+2} - \theta_j) \\
\nonumber & + & \frac{1}{2} \left( 4 + 2
\cos(\theta_{j+1}-\theta_{j+2}) + 2 \cos(\theta_{j-1} -
\theta_{j-2}) + \cos(\theta_j - \theta_{j-2})
\right) \sin(\theta_{j-2} - \theta_j) \\
\nonumber & + & 2 ( 1 + \cos(\theta_{j+1} - \theta_j)) \sin(\theta_j
- \theta_{j+1}) + 2 ( 1 + \cos(\theta_{j-1} - \theta_j))
\sin(\theta_j - \theta_{j-1}) \\
\nonumber & + & 2 (2 + \cos(\theta_2 - \theta_1) + \cos(\theta_3 -
\theta_1) +\cos(\theta_4 - \theta_1) +\cos(\theta_3 - \theta_2)
+\cos(\theta_4 - \theta_2) +\cos(\theta_4 - \theta_3) ) \\
& \phantom{t} & \times ( \sin(\theta_j - \theta_{j+1}) +
\sin(\theta_j - \theta_{j-1}) + \sin(\theta_j - \theta_{j+2})) \\
\nonumber & + & 4 \sin(\theta_j - \theta_{j+1}) + 4 \sin(\theta_j -
\theta_{j-1})
\end{eqnarray*}
} For the asymmetric vortex, we have
$$
g^{(4)}_j = (-1)^j 2 \sin(2 \theta), \qquad j = 1,2,3,4.
$$
The Jacobian matrix ${\cal M}_2$ has two zero eigenvalues with
orthogonal eigenvectors:
$$
{\bf p}_1 = \left( \begin{array}{c} 1 \\ 1 \\ 1 \\ 1 \end{array}
\right), \qquad {\bf p}_2 = \left( \begin{array}{c} -1 \\ 1 \\ -1 \\
1 \end{array} \right).
$$
It is clear that the vector ${\bf g}^{(4)} = 2 \sin(2 \theta) {\bf
p}_2$ is not orthogonal to the eigenvector ${\bf p}_2$ of the
kernel of ${\cal M}_2$, unless $\theta = \{0,\frac{\pi}{2},\pi\}$.
By Proposition 2.10 in \cite{paperII}, the family of asymmetric
vortices (\ref{asymmetric-vortices}) terminates at the
fourth-order reduction. The exceptional cases include discrete
solitons for $\theta = \{ 0, \pi\}$ and the vortex cross at
$\theta = \frac{\pi}{2}$. In order to consider persistence of the
vortex cross, we compute the Jacobian matrices ${\cal M}_2$ and
${\cal M}_4$ from the bifurcation functions ${\bf g}^{(2)}$ and
${\bf g}^{(4)}$ explicitly,
$$
{\cal M}_2 = \left( \begin{array}{ccccc} -1 & 0 & 1 & 0
\\ 0 & -1 & 0 & 1 \\ 1 & 0 & -1 & 0 \\ 0 & 1 & 0 & -1
\end{array} \right), \qquad {\cal M}_4 = \left( \begin{array}{ccccc} 3
& 2 & -7 & 2 \\ 2 & 3 & 2 & -7 \\
-7 & 2 & 3 & 2 \\ 2 & -7 & 2 & 3 \end{array} \right)
$$
Since ${\cal M}_4 {\bf p}_1 = {\bf 0}$ and ${\cal M}_4 {\bf p}_2
\neq {\bf 0}$, the zero eigenvalue of ${\cal M}_2$ with the
associated eigenvector ${\bf p}_2$ bifurcates. By Proposition 2.9
in \cite{paperII}, this implies that the family of the vortex
cross is continued from the anti-continuum limit uniquely up to
the rotational transformation $\mbox{\boldmath $\theta$} \to
\mbox{\boldmath $\theta$} + \theta_0 {\bf p}_1$ that corresponds
to the gauge symmetry of the dNLS equation (\ref{scalar-eq}).
Proposition \ref{proposition-existence} is hence proved.

Small eigenvalues of the linearized Jacobian matrix ${\cal
H}(\epsilon)$ are defined by an extended eigenvalue problem for
the Jacobian matrices ${\cal M}_2$ and ${\cal M}_4$,
$$
\left( \epsilon^2 {\cal M}_2 + \epsilon^4 {\cal M}_4 + {\rm
O}(\epsilon^6) \right) {\bf c} = \gamma {\bf c}.
$$
There exist four eigenvalues of the extended problem which admit
the asymptotic approximations,
$$
\gamma_{1,2} = -2 \epsilon^2 + {\rm O}(\epsilon^4),  \qquad \gamma_3
= -8 \epsilon^4 + {\rm O}(\epsilon^6), \qquad   \gamma_4 = 0.
$$
The eigenvalue $\gamma_3$ is obtained by the perturbation theory for
the zero eigenvalue of ${\cal M}_2$ associated with the eigenvector
${\bf p}_2$ (orthogonal to the eigenvector ${\bf p}_1$):
$$
\lim_{\epsilon \to 0} \epsilon^{-4} \gamma_3 = \frac{({\bf p}_2,
{\cal M}_4 {\bf p}_2)}{({\bf p}_2,{\bf p}_2)} = -8.
$$

\end{document}